\documentclass[review]{elsarticle}
\usepackage{amsmath,amssymb,amsfonts,xcolor}
\usepackage {longtable}

\usepackage{lineno,hyperref}
\modulolinenumbers[5]

\journal{J. Comp. \&  Appl. Math., ISSN: 0377-0427}
\newtheorem{theorem}{Theorem}%[section]
\newtheorem{lemma}{Lemma}%[section]
\newtheorem{remark}{Remark}%[section]
\newtheorem{example}{Example}%[section]

\begin{document}

\begin{frontmatter}

\title{Exponentially convergent symbolic algorithm of the functional-discrete method for the fourth order Sturm--Liouville problems with polynomial coefficients}

\author{Volodymyr Makarov\fnref{myfootnote}}
\author{Nataliia Romaniuk\corref{mycorrespondingauthor}}
\address{Department of Numerical Mathematics, Institute of Mathematics of National Academy of Sciences of Ukraine, 3 Tereshchenkivs'ka Str., 01004 Kyiv-4, Ukraine}
\fntext[myfootnote1]{{\textit{Email address:}} makarov@imath.kiev.ua (Volodymyr Makarov)}

\cortext[mycorrespondingauthor]{Corresponding author}
\ead{romaniuknm@gmail.com}

\begin{abstract}
A new symbolic algorithmic implementation of the functional-discrete (FD-) method is developed and justified for the solution of fourth order Sturm--Liouville problem on a finite interval in the Hilbert space. The eigenvalue problem for the fourth order ordinary differential equation with polynomial coefficients is investigated. The sufficient conditions of an exponential convergence rate of the proposed approach are received. The obtained estimates of the absolute errors of FD-method significantly improve the accuracy of the estimates obtained earlier by I.P~Gavrilyuk, V.L.~Makarov and A.M.~Popov in 2010. Our algorithm is symbolic and operates with the decomposition coefficients of the eigenfunction corrections in some basis. The number of summands in these decompositions depends on the degree of the potential coefficients and the correction number. Our method uses only the algebraic operations and basic operations on $(2\times 1)$ column vectors and $(2\times 2)$ matrices. The proposed approach does not require solving any boundary value problems and computations of any integrals, unlike the previous variants of FD-method by I.P.~Gavrilyuk, V.L.~Makarov, A.M.~Popov and N.M.~Romaniuk in 2010 and 2017. The corrections to eigenpairs are computed exactly as analytical expressions, and there are no rounding errors. The numerical examples illustrate the theoretical results. The numerical results obtained with the FD-method are compared with the numerical test results obtained with other existing numerical techniques.

\end{abstract}


\end{frontmatter}

%\linenumbers

\section{Introduction}\label{MakRom_section_1}

There is a great number of numerical methods for Sturm--Liouville problems for the second- and higher-order ordinary differential equations. The analytical methods based on perturbation and homotopy ideas \cite{115_GMR_Arm1979,118_GMR_AllgGeorg1990} are widely used for solving the eigenvalue problems.
The numerical-analytical (functional-discrete) methods refer to these methods (see, for example, \cite{Mak1991,
MakarovKlymenko2007,
GMR-GavrMakPop2010,
MakRom2014Dir,
Conf_MakRom2015,
MakRom2017,
GavrMakRom2017,
GavrMakRom20152017} and compare with Adomian decomposition method  \cite{MRL-Adomian1993,65_Rach2012}).  
Using these analytical methods the solutions can be found as fast convergent functional series. The properties of the solution of the original problem can be investigated with the help of approximation solution. Moreover, these approaches allow in a natural way to use computer algebra systems for developing and implementation of symbolic-numerical algorithms.

The functional-discrete (FD-) method was suggested by V.~Makarov \cite{Mak1991} in 1991. In subsequent years, FD-method was developed for the solution of many different problems. This method enables us to overcome many disadvantages of the discrete methods such as: 
\begin{enumerate}
\item[--] the accuracy degradation with the increasing of the eigenvalue index;
\item[--] usage of the mesh generated at the start of the numerical process;
\item[--] saturation of accuracy;
\item[--] the number of reliable numerical eigenvalues is limited and depends on a mesh step (see \cite{Pryce1993,Zhang2015}).
\end{enumerate}
The main advantages of FD-method are its features, which differ from many other methods: 
\begin{enumerate}
\item the approach can be applied to operator equations in general form; 
\item the approach can be applied also to eigenvalue problems with multiple eigenvalues (see, for example, \cite{GavrMakRom20152017}); 
\item all eigenpairs can be computed in parallel; 
\item the convergence rate increases as the index of the eigenpair increases; 
\item it was proved that in many cases the FD-method converges exponentially or super-exponentially.
\end{enumerate}
 
Presented modification of the traditional algorithm of the functional-discrete (FD-) method was proposed in \cite{MakarovKlymenko2007}. The general idea of the symbolic algorithms for FD-method is the representation of the eigenfunction corrections in some basis. Then we obtain recurrence relations for the decomposition coefficients. 
Finally, our algorithm operates with these decomposition coefficients. The modified FD-method does not require solving any boundary value problems and computations of any integrals. In certain cases, the algorithm uses only the algebraic operations. 
Moreover, the corrections to eigenpairs are computed exactly as analytical expressions, and there are no rounding errors.
But if computational difficulties or memory overflow arise then we can avoid combinatorial explosion. In proposed approach instead of using rational arithmetic we can easily transit to floating--point arithmetic which "represents an alternative idea: round the computation at every step, not just at the end" (see \cite{Trefethen2015}).

Briefly described general idea was used for developing and justification of new symbolic algorithms for FD-method for the Sturm--Liouville problems on a finite interval for the Schr{\"o}dinger equation with a polynomial potential in \cite{MakRom2014Dir,MakRom2017}. 
Using the described idea for the fourth order Sturm--Liouville problem, in this article we modify the traditional method from \cite{GMR-GavrMakPop2010,GavrMakRom20152017} and develop a new symbolic algorithm of the FD-method.
The proposed algorithm of our method is developed when the potential coefficients are approximated by zero function. For this case, FD-method is purely analytical method and may be considered one of the variants of the homotopy method \cite{115_GMR_Arm1979,118_GMR_AllgGeorg1990}. 
Note that some of the results of this article were announced in \cite{Conf_MakRom2015}. Unlike the symbolic algorithm in \cite{Conf_MakRom2015}, the presented approach produces explicit recursive formulas for the decomposition coefficients of the representation for the eigenfunctions corrections.

The article is organized as follows. Section~\ref{MakRom_section_2} deals with the problem statement. 
Section~\ref{MakRom_section_3} contains the traditional algorithm of the simplest variant of the FD-method. 
In Section~\ref{MakRom_section_4} a new structural representation of the eigenfunctions corrections is obtained. This representation is used in Section~\ref{MakRom_section_6} to develop a new symbolic algorithm of the FD-method. 
The sufficient conditions of an exponential convergence rate and the absolute errors estimates of the proposed approach are received in Section~\ref{MakRom_section_5}. The obtained absolute errors estimates of the FD-method significantly improve the accuracy of the estimates obtained earlier in \cite{GMR-GavrMakPop2010}.
Derivation of the basic formulas for the proposed new symbolic algorithmic implementation of our method is given in Section~\ref{MakRom_section_6} and in Appendices A, B, C. 
The numerical algorithm is given in Section~\ref{MakRom_section_6_1}.
Section~\ref{MakRom_section_7} illustrates the theoretical results by numerical examples. 
In Examples~\ref{example_MakRom_1} and \ref{example_MakRom_2} the numerical results obtained with the FD-method are compared with the numerical test results obtained with other existing numerical techniques \cite{58_AttiliLesnic2006,25_SyamSiyyam2009,HPM2010Ex1,HAM2011Ex1,59_Chanane2010,RATTANA2013144}.
A review of the obtained results with implementation features and advantages of our developed numerical method are given in the final Section~\ref{MakRom_section_8}.

\section{Problem statement}\label{MakRom_section_2}%2

Let us consider the regular Sturm--Liouville problem in a Hilbert space for the fourth order ordinary differential equation 
\begin{equation} \label{GrindEQ__1_} 
u^{{\rm IV}} (x)+q_{2} (x)u''(x)+q_{1} (x)u'(x)+\left(q_{0} (x)-\lambda \right)u(x)=0,
\end{equation} 
$$
x\in {\kern 1pt} (0,X),\;\;0<X<\infty , 
$$
with the boundary conditions
\begin{equation} \label{GrindEQ__2_} u(0)=u''(0)=u(X)=u''(X)=0, 
\end{equation} 
where $X$ is the real constant. The real-valued polynomial coefficients are
\begin{equation} \label{GrindEQ__3_} 
q_{0} (x)=\sum _{l=0}^{r_{0} }A_{l}  x^{l} ,\; \;
q_{1} (x)=\sum _{l=0}^{r_{1} }B_{l}  x^{l} ,\; \;  
q_{2} (x)=\sum _{l=0}^{r_{2} }C_{l}  x^{l} ,\; \; r_{i} \ge 1, i=0,1,2, 
\end{equation} 
where the constants $r_{i}$, $i=0,1,2$ are positive integers, and $q_{i} (x)\in C^{i} [0,X]$, $i=0,1,2$.

In \cite{GMR-GavrMakPop2010} it was shown that the fourth order ordinary differential equation with all derivatives of the eigenfunction could be reduced to the form (\ref{GrindEQ__1_}) using the variable transformation.
 That is why we consider the eigenvalue problem with equation (\ref{GrindEQ__1_}).

\section{Traditional algorithm of the simplest variant of the FD-method}\label{MakRom_section_3}%3

In this article, the simplest variant of the FD-method is applied to the Sturm--Liouville problem (\ref{GrindEQ__1_})--(\ref{GrindEQ__3_}). It means that we consider the simplest case of the approximation of potential coefficients (\ref{GrindEQ__3_}) by zero function
\[\bar{q}_{s} (x)\equiv 0,\;\;\;s=0,1,2.\]
This section contains the traditional algorithm of the FD-method which in this case is the purely analytical method. In this case, FD-method can be considered as one of the variants of the homotopy method (see \cite{115_GMR_Arm1979,118_GMR_AllgGeorg1990}), and its idea is closely related to the ideas of the Adomian decomposition method \cite{MRL-Adomian1993,65_Rach2012}. Below the symbolic algorithm for the simplest variant of the FD-method is developed.

Note that in the case when the simplest variant of the FD-method (with $\bar{q}_{s} (x)\equiv 0$, $s=0,1,2$)  is divergent for the smallest eigenvalues of the problem (\ref{GrindEQ__1_})--(\ref{GrindEQ__3_}), the general scheme of the FD-method (usually) with piecewise-constant approximations to potential coefficients is used. Developing and justification of a symbolic algorithmic implementation of the general scheme of the FD-method for the problem (\ref{GrindEQ__1_})--(\ref{GrindEQ__3_}) are slated for the near future.

The exact solution of the eigenvalue problem (\ref{GrindEQ__1_})--(\ref{GrindEQ__3_}) is then represented by the series
\begin{equation} \label{GrindEQ__3_11} 
u_{n} (x)=\sum _{j=0}^{\infty } u_{n}^{(j)} (x),{\kern 1pt} \; \; \; \lambda _{n} =\sum _{j=0}^{\infty } \lambda _{n}^{(j)} ,
\end{equation} 
provided that these series converge. The sufficient conditions for the convergence of the series (\ref{GrindEQ__3_11}) will be presented later in Section~\ref{MakRom_section_5}. The approximate solution to the problem (\ref{GrindEQ__1_})--(\ref{GrindEQ__3_}) is represented by a pair of corresponding truncated series, namely,
\begin{equation} \label{GrindEQ__4_} 
\mathop{u _{n}}\limits^{m} (x)=\sum _{j=0}^{m} \; u_{n}^{(j)} (x),{\kern 1pt} \; \; \; \mathop{\lambda_{n} }\limits^{m}  =\sum _{j=0}^{m} \; \lambda _{n}^{(j)} , 
\end{equation} 
which is called an approximation of rank $m$ to eigenpair (eigenfunction $u_{n}(x)$ and eigenvalue $\lambda _{n}$) with index number $n$. 
The summands of series (\ref{GrindEQ__4_}) $u_{n}^{(j)}  (x)$, $\lambda _{n}^{(j)}$ are called the corrections to eigenpairs at the $j$--th step of the FD-method. The corrections $u_{n}^{(j+1)}  (x)$, $\lambda _{n}^{(j+1)}$ are the solutions of the following recursive sequence of problems (see \cite{GMR-GavrMakPop2010})
\begin{equation} \label{GrindEQ__5_} 
\frac{d^{4} }{dx^{4} } u_{n}^{(j+1)} (x)-\lambda _{n}^{(0)} u_{n}^{(j+1)} (x)=F_{n}^{(j+1)} (x),\; \; \; 
x\in (0,X), 
\end{equation} 
\begin{equation} \label{GrindEQ__6_} 
u_{n}^{(j+1)} (0)=\frac{d^{2} u_{n}^{(j+1)} (0)}{dx} =u_{n}^{(j+1)} (X)=\frac{d^{2} u_{n}^{(j+1)} (X)}{dx} =0,\; \end{equation} 
$$j=0,1,...,m-1,$$
where
\begin{equation} \label{GrindEQ__6_1}
\begin{split}
F_{n}^{(j+1)} (x)&=\sum _{p=0}^{j}\lambda _{n}^{(j+1-p)} u_{n}^{(p)} (x)-q_{2} (x)\frac{d^{2} }{dx^{2} } u_{n}^{(j)} (x)\\
&-q_{1} (x)\frac{d}{dx} u_{n}^{(j)} (x)- q_{0} (x)u_{n}^{(j)} (x).
\end{split}
\end{equation} 
Using the solvability condition 
\begin{equation}\label{GrindEQ__6_2}
\left(F_{n}^{(j+1)} (x),u_{n}^{(0)} (x)\right)_{L_{2} \left(0,X\right)} =\int _{0}^{X} F_{n}^{(j+1)} (x)u_{n}^{(0)} (x) dx=0
\end{equation}
of problems (\ref{GrindEQ__5_})--(\ref{GrindEQ__6_1}) for a fixed $j$ $(j = 0, 1, . . .)$, we obtain the following formula for the eigenvalue corrections:
\begin{equation} \label{GrindEQ__7_} 
\lambda _{n}^{(j+1)} =\int _{0}^{X}\left(q_{2} (x)\frac{d^{2} }{dx^{2} } u_{n}^{(j)} (x)+q_{1} (x)\frac{d}{dx} u_{n}^{(j)} (x)+q_{0} (x)u_{n}^{(j)} (x)\right)u_{n}^{(0)} (x) dx. 
\end{equation} 
Solutions to the problems (\ref{GrindEQ__5_})--(\ref{GrindEQ__6_1}) satisfy the orthogonality condition
\begin{equation} \label{GrindEQ__8_} 
\left(u_{n}^{(j+1)} (x),u_{n}^{(0)} (x)\right)_{L_{2} \left(0,X\right)} =\int _{0}^{X}  u_{n}^{(j+1)} (x)u_{n}^{(0)} (x) dx=0,  
\end{equation} 
$$j=0,1,...,m-1.$$
The initial approximation $u_{n}^{(0)} (x)$, $\lambda _{n}^{(0)}$ is the solution of the so-called base problem, that is,
\begin{equation} \label{GrindEQ__9_} 
\frac{d^{4} u_{n}^{(0)} (x)}{dx^{4} } -\lambda _{n}^{(0)} u_{n}^{(0)} (x)=0,\;\;x\in (0,X), 
\end{equation} 
\begin{equation} \label{GrindEQ__9_1} 
u_{n}^{(0)} (0)=\frac{d^{2} u_{n}^{(0)} (0)}{dx} =u_{n}^{(0)} (X)=\frac{d^{2} u_{n}^{(0)} (X)}{dx} =0.
\end{equation} 
The solution of (\ref{GrindEQ__9_})--(\ref{GrindEQ__9_1}) is the following
\begin{equation} \label{GrindEQ__10_} 
u_{n}^{(0)} (x)=\sqrt{\frac{2}{X} } \sin \left(\frac{n\pi }{X} x\right), \;\;\lambda _{n}^{(0)} =\frac{\left(n\pi \right)^{4} }{X^{4} }.  
\end{equation}

\section{Representation of the corrections to eigenfunctions $u_{n}^{(j)}  (x)$}\label{MakRom_section_4}%4

Let us introduce the generalized Green's function for a linear differential operator, 
corresponding to the problems by (\ref{GrindEQ__5_})--(\ref{GrindEQ__10_}), in the following form
\begin{equation} \label{GrindEQ__11_} 
g_{n} (x,\xi )=\frac{2X^{3} }{\pi ^{4} } \sum_{p=1,p\ne n}^{\infty } \frac{\sin \left(\frac{p\pi }{X} x\right)\sin \left(\frac{p\pi }{X} \xi \right)}{p^{4} -n^{4} } ,\; \; x,\xi \in \left[0,X\right]. 
\end{equation} 
The function (\ref{GrindEQ__11_}) can be expressed as
\begin{equation} \label{GrindEQ__12_}
\begin{split}
g_{n} (x,\xi )&=\frac{X^{3} }{2\pi ^{3} n^{3} } \left[-\left(\frac{x}{X} -{\rm H} \left(\frac{x-\xi }{X} \right)\right)\cos \left(\frac{n\pi }{X} x\right)\right. \sin \left(\frac{n\pi }{X} \xi \right)\\
&-\left(\frac{\xi }{X} -{\rm H} \left(\frac{\xi -x}{X} \right)\right)\sin \left(\frac{n\pi }{X} x\right)\cos \left(\frac{n\pi }{X} \xi \right)+\\
&+\frac{1}{\sinh \left(\pi n\right)} \sinh \left(\pi n\left(\frac{x}{X} -{\rm H} \left(\frac{x-\xi }{X} \right)\right)\right)\\
&\times\sinh \left(\pi n\left(\frac{\xi }{X} -{\rm H} \left(\frac{\xi -x}{X} \right)\right)\right)\left. +\frac{3}{2\pi n} \sin \left(\frac{n\pi }{X} x\right)\sin \left(\frac{n\pi }{X} \xi \right)\right], 
\end{split}
\end{equation} 
$$x,\xi \in \left[0,X\right],$$
where ${\rm H}(x)$ is the Heaviside function, and ${\rm H} (0)=1$.

\begin{lemma}\label{MakRom_lemma_1} 
The generalized Green's function (\ref{GrindEQ__11_}), (\ref{GrindEQ__12_}) has the following properties:
$$
g_{n} (x,\xi )=g_{n} (\xi ,x), 
\;\;\;\;\;
g_{n} (x,\xi )=g_{n} (X-x,X-\xi ),
$$
$$
\int _{0}^{X} g_{n} (x,\xi )\sin \left(\frac{\pi n}{X} x\right)dx=0. 
$$
\end{lemma}
For a fixed $j$ the solution of a problem (\ref{GrindEQ__5_})--(\ref{GrindEQ__6_2}), (\ref{GrindEQ__3_}), which satisfies the orthogonality condition (\ref{GrindEQ__8_}), can be expressed by a formula
\begin{equation} \label{GrindEQ__14_} 
u_{n}^{(j+1)} (x)=\int _{0}^{X} g_{n} (x,\xi )F_{n}^{(j+1)} (\xi ) d\xi . 
\end{equation} 

The following assertion is proved by a method of complete induction.
To accomplish this we use the integral representation (\ref{GrindEQ__14_}) and the solution of the base problem (\ref{GrindEQ__10_}),
as well as the properties of the problems (\ref{GrindEQ__5_})--(\ref{GrindEQ__6_2}), (\ref{GrindEQ__3_}) and of the generalized Green's function (\ref{GrindEQ__11_}), (\ref{GrindEQ__12_}) (see Lemma~\ref{MakRom_lemma_1}).
\begin{lemma}\label{MakRom_lemma_2} 
The solution of problem (\ref{GrindEQ__5_})--(\ref{GrindEQ__6_2}), (\ref{GrindEQ__3_}) can be represented by
\begin{equation} \label{GrindEQ__15_} 
\begin{split}
u_{n}^{(j+1)}(x)&=\sum _{p=0}^{M(j+1)} x^{p} \left[b_{n,p}^{(j+1)} \cos \left(\frac{\pi n}{X} x\right)+a_{n,p}^{(j+1)} \sin \left(\frac{\pi n}{X} x\right)\right]\\
&+\sum _{p=0}^{M(j)} x^{p} \left[d_{n,p}^{(j+1)} \cosh \left(\frac{\pi n}{X} x\right)+c_{n,p}^{(j+1)} \sinh \left(\frac{\pi n}{X} x\right)\right],\\
&j=0,1,...,\;n=1,2,...
\end{split}
\end{equation} 
where $M(j)=j(r+1)$, $r=\max \left\{r_{0} ,r_{1} ,r_{2} \right\},$ $a_{n,0}^{(0)}=\sqrt{\frac{2}{X} } $, $b_{n,0}^{(0)} =d_{n,0}^{(0)} =c_{n,0}^{(0)} =0$. 
 \end{lemma}
 
In Lemma~\ref{MakRom_lemma_2} the coefficients 
\begin{equation}\label{GrindEQ__66_}  
\begin{split}
&b_{p}^{(j+1)} ,\;\;a_{p}^{(j+1)} \;\;(p=0,1,...,M(j+1)),\\
&c_{s}^{(j+1)} ,\;\;d_{s}^{(j+1)}\;\;(s=0,1,...,M(j)) 
\end{split}
\end{equation}
are the decomposition coefficients of the eigenfunction corrections $u_{n}^{(j+1)} (x)$ in the basis
$$x^{p} \cos \left(\frac{\pi n}{X} x\right),\;\;x^{p} \sin \left(\frac{\pi n}{X} x\right)\;\;(p=0,1,...,M(j+1)),$$
$$x^{s} \cosh \left(\frac{\pi n}{X} x\right),\;\;x^{s} \sinh \left(\frac{\pi n}{X} x\right)\;\;(s=0,1,...,M(j))$$
on interval $\left[0,X\right]$. 
Unlike \eqref{GrindEQ__14_}, the representation  \eqref{GrindEQ__15_} is used below to develop a new symbolic algorithmic implementation of the FD-method, numerical algorithm for which is fully given in Section~\ref{MakRom_section_6_1}. Exact explicit recursive formulas for these coefficients are found in Section~\ref{MakRom_section_6} (see \eqref{GrindEQ__65_1_}, \eqref{GrindEQ__65_2_}, \eqref{GrindEQ__66_1_} and \eqref{GrindEQ__66_2_}).

\section{Convergence of the FD-method}\label{MakRom_section_5}%5

To investigate the convergence of the FD-method we substitute the expressions (\ref{GrindEQ__7_}) and (\ref{GrindEQ__14_}) into formulas
\begin{equation} \label{GrindEQ__16_} 
\begin{split}
\lambda _{n}^{(j+1)}& =\sqrt{\frac{2}{X} } \int _{0}^{X} \left(-\left(\frac{n\pi }{X} \right)^{2} \right. q_{2} (x)\sin \left(\frac{n\pi }{X} x\right)\\
&+\frac{n\pi }{X} \left[2q'_{2} (x)-q_{1} (x)\right]\cos \left(\frac{n\pi }{X} x\right)\\
&+\left[q''_{2} (x)-q'_{1} (x)+q_{0} (x)\right]\left. \sin \left(\frac{n\pi }{X} x\right)\right)u_{n}^{(j)} (x)dx
\end{split}
\end{equation} 
and
\begin{equation} \label{GrindEQ__17_} 
\begin{split}
u_{n}^{(j+1)} &(x)=\frac{2X^{3} }{\pi ^{4} } \sum_{p=1,p\ne n}^{\infty } \frac{1}{p^{4} -n^{4} } \sin \left(\frac{p\pi }{X} x\right)\\
&\times\int _{0}^{X} \left\{\left[\left(\frac{p\pi }{X} \right)^{2} q_{2} (\xi )\sin \left(\frac{p\pi }{X} \xi \right)\right. \right.\\
&+\frac{p\pi }{X} \left[-2q'_{2} (\xi )+q_{1} (\xi )\right]\cos \left(\frac{p\pi }{X} \xi \right)\\
&+\left[-q''_{2} (\xi )+q'_{1} (\xi )-q_{0} (\xi )\right]\left. \sin \left(\frac{p\pi }{X} \xi \right)\right]u_{n}^{(j)} (\xi )\\
&+\left. \sum _{s=0}^{j} \; \lambda _{n}^{(j+1-s)} u_{n}^{(s)} (\xi )\sin \left(\frac{p\pi }{X} \xi \right)\right\}d\xi .
\end{split}
\end{equation} 
The formulas (\ref{GrindEQ__16_}) and (\ref{GrindEQ__17_}) were obtained using the integration by parts as well as the representation for the generalized Green's function $g_{n} (x,\xi )$ by a series (\ref{GrindEQ__11_}). Note that obtained in this Section absolute errors estimates of the FD-method significantly improve the accuracy of the estimates obtained earlier in \cite{GMR-GavrMakPop2010}.

One can deduce from (\ref{GrindEQ__16_}) the next estimate for the eigenvalue corrections
\begin{equation} \label{GrindEQ__18_} 
\left|\lambda _{n}^{(j+1)} \right|\le \omega \sqrt{\frac{2}{X} } \left[\left(\frac{n\pi }{X} \right)^{2} +\frac{n\pi }{X} +1\right]\left\| u_{n}^{(j)} \right\| , 
\end{equation} 
where 
\[\omega =\max \left\{\left\| q_{2} \right\| _{\infty } ,\left\| 2q'_{2} -q_{1} \right\| _{\infty } ,\left\| q''_{2} -q'_{1} +q_{0} \right\| _{\infty } \right\},\;\;\;\;\left\| v\right\| _{\infty } =\mathop{\max }\limits_{x\in [0,X]} \left|v(x)\right|,\] 
\[\left\| v\right\| =\sqrt{\left(v(x),v(x)\right)_{L_{2} \left(0,X\right)} } =\left(\int _{0}^{X} \; \left[v(x)\right]^{2} dx\right)^{^{{1\mathord{\left/ {\vphantom {1 2}} \right. \kern-\nulldelimiterspace} 2} } } .\]  
It is easy to establish that the following inequalities are correct:
\[\sqrt{\sum _{p=1,p\ne n}^{\infty } \frac{p^{4-2k} }{\left(p^{4} -n^{4} \right)^{2} } } \le \frac{n^{1-k} }{2n^{2} -2n+1} ,\, k=0,1,2,\, n=1,2,...\] 
which are used to obtain from (\ref{GrindEQ__18_}) the estimate for the eigenfunction corrections (\ref{GrindEQ__17_}):    
\begin{equation} \label{GrindEQ__19_} 
\begin{split}
&\left\| u_{n}^{(j+1)} \right\| \le \frac{X^{2} }{\pi ^{2} } \frac{1}{2n^{2} -2n+1} \left[\left[n\left\| q_{2} \right\| _{\infty } +\frac{X}{\pi } \left\| -2q'_{2} +q_{1} \right\| _{\infty }\right. \right.\\
&\left. +\frac{X^{2} }{\pi ^{2} } \frac{\left\| -q''_{2} +q'_{1} -q_{0} \right\| _{\infty } }{n} \right]\left\| u_{n}^{(j)} \right\|+\sqrt{\frac{2}{X} } \omega  \left[n+\frac{X}{\pi } +\frac{X^{2} }{n\pi ^{2} } \right]\\
&\times\left.\sum _{s=0}^{j} \; \left\| u_{n}^{(j-s)} \right\| \left\| u_{n}^{(s)} \right\| \right]\le \frac{X^{2} }{\pi ^{2} } \frac{\omega }{2n^{2} -2n+1} \left[n+\frac{X}{\pi } +\frac{X^{2} }{n\pi ^{2} } \right]\\
&\times\left[\left\| u_{n}^{(j)} \right\| +\sqrt{\frac{2}{X} } \sum _{s=0}^{j} \; \left\| u_{n}^{(j-s)} \right\| \left\| u_{n}^{(s)} \right\| \right]\le M_{n} \sum _{s=0}^{j} \; \left\| u_{n}^{(j-s)} \right\| \left\| u_{n}^{(s)} \right\| ,
\end{split}
\end{equation} 
where
\begin{equation} \label{GrindEQ__20_} M_{n} =\frac{X^{2} }{\pi ^{2} } \frac{\omega }{2n^{2} -2n+1} \left[n+\frac{X}{\pi } +\frac{X^{2} }{n\pi ^{2} } \right]\max \left\{1,\sqrt{\frac{2}{X} } \right\}. \end{equation} 
Substituting in (\ref{GrindEQ__19_}) 
 \begin{equation} \label{GrindEQ__21_} 
U_{j} =M_{n}^{-j} \left\| u_{n}^{\left(j\right)} \right\| ,\;\;\;\;U_{0} =\left\| u_{n}^{\left(0\right)} \right\| =1 
\end{equation} 
and replacing the new variables by the majorant variables subject to 
$$U_{j} \le \overline{U}_{j},\;\;\;\;\overline{U}_{0} =U_{0} =1,$$
we come to the majorant equation 
\begin{equation} \label{GrindEQ__22_} 
\overline{U}_{j+1} =\sum _{s=0}^{j}\overline{U}_{j-s} \overline{U}_{s}. 
\end{equation} 
This equation is the nonlinear recurrence relation and the so-called convolution-type equation.
The solution of the equation (\ref{GrindEQ__22_}) is (see, e.g., \cite[p.~159-161,210]{Vilenkin1969}, \cite{ReinNievDeo1977})
\begin{equation} \label{GrindEQ__23_} 
\overline{U}_{j+1} =\frac{\left(2j+2\right)!}{\left(j+1\right)!\left(j+2\right)!} =4^{j+1} 2\frac{\left(2j+1\right)!!}{\left(2j+4\right)!!} . 
\end{equation} 
Returning to the old variables (see substitution of variables (\ref{GrindEQ__21_})), we obtain from (\ref{GrindEQ__23_}) the following estimate for the solution of (\ref{GrindEQ__19_}):
\begin{equation} \label{GrindEQ__24_} 
 \left\| u_{n}^{\left(j+1\right)} \right\| \le \left(4M_{n} \right)^{j+1} 2\frac{\left(2j+1\right)!!}{\left(2j+4\right)!!} \le \frac{\left(4M_{n} \right)^{j+1} }{\left(j+2\right)\sqrt{\pi \left(j+1\right)} },
\end{equation} 
and then, from (\ref{GrindEQ__18_}), the next estimate for the eigenvalue corrections
\begin{equation} \label{GrindEQ__25_} 
\begin{split}
\left|\lambda _{n}^{(j+1)} \right|&\le \omega \sqrt{\frac{2}{X} } \left[\left(\frac{n\pi }{X} \right)^{2} +\frac{n\pi }{X} +1\right]\left(4M_{n} \right)^{j} 2\frac{\left(2j-1\right)!!}{\left(2j+2\right)!!}\\
&\le \omega \sqrt{\frac{2}{X} } \left[\left(\frac{n\pi }{X} \right)^{2} +\frac{n\pi }{X} +1\right]\frac{\left(4M_{n} \right)^{j} }{\left(j+1\right)\sqrt{\pi j} } . 
\end{split}
\end{equation} 
The last parts of inequalities (\ref{GrindEQ__24_}) and (\ref{GrindEQ__25_}) were obtained using the judgements like those from the proof of the Wallis formula (see, e.g., \cite[p.~344]{Fichtenholz1968}).
From estimates (\ref{GrindEQ__24_}) and (\ref{GrindEQ__25_}) follows the next theorem which contains the sufficient conditions of an exponential convergence rate of the FD-method and estimates of its absolute errors.
\begin{theorem}\label{MakRom_theorem_3}
Let $q_{i} (x)\in C^{i} [0,X]$, $i=0,1,2$ and let the following condition hold true:
\begin{equation} \label{GrindEQ__26_} 
r_{n} =4M_{n} <1,  n=1,2,... .                                                
\end{equation} 
Then the FD-method for the Sturm--Liouville problem (\ref{GrindEQ__1_})--(\ref{GrindEQ__3_}) converges exponentially and the following estimates of the absolute errors are valid:
\begin{equation} \label{GrindEQ__27_} 
\left|\lambda _{n} -\mathop{\lambda_{n} }\limits^{m}  \right|\le \omega \sqrt{\frac{2}{X} } \left[\left(\frac{n\pi }{X} \right)^{2} +\frac{n\pi }{X} +1\right]\frac{\left(r_{n} \right)^{m} }{1-r_{n} } \frac{1}{\left(m+1\right)\sqrt{\pi m} } , 
\end{equation} 
\begin{equation} \label{GrindEQ__28_} 
\left\| u_{n} -\mathop{u_{n}}\limits^{m}  \right\| \le 2\frac{\left(r_{n} \right)^{m+1} }{1-r_{n} } \frac{\left(2m+1\right)!!}{\left(2m+4\right)!!} \le \frac{\left(r_{n} \right)^{m+1} }{\left(m+2\right)\sqrt{\pi \left(m+1\right)} } . 
\end{equation} 
\end{theorem}

\section{Derivation of basic formulas for the symbolic algorithm of the FD-method}\label{MakRom_section_6}%6

Further, in this section, we describe and develop a new symbolic algorithm of the FD-method for the problem (\ref{GrindEQ__1_})--(\ref{GrindEQ__3_}). In this section, basic formulas of the symbolic algorithm are given and also additional formulas are given in Appendices~A, B, C.
Unlike the symbolic algorithm from \cite{Conf_MakRom2015}, the presented approach produces explicit recursive formulas for the coefficients in (\ref{GrindEQ__15_}) at the $(j+1)$--th step of the FD-method. Then the computer algebra system Maple was used for a software implementation.

Let us substitute (\ref{GrindEQ__15_}) into (\ref{GrindEQ__6_1}), (\ref{GrindEQ__3_}) and group together the summands as follows 
\begin{equation} \label{GrindEQ__29_1}
\begin{split}
F_{n}^{(j+1)} (x)&=F_{n,\cos }^{(j+1)} (x)\cos \left(\frac{\pi n}{X} x\right)+F_{n,\sin }^{(j+1)} (x)\sin \left(\frac{\pi n}{X} x\right)\\
&+F_{n,\cosh }^{(j+1)} (x)\cosh \left(\frac{\pi n}{X} x\right)+F_{n,\sinh }^{(j+1)} (x)\sinh \left(\frac{\pi n}{X} x\right)
\end{split}
\end{equation} 
changing the order of summation in analytical expressions for $F_{n,\cos }^{(j+1)} (x)$, $F_{n,\sin }^{(j+1)} (x)$, $F_{n,\cosh }^{(j+1)} (x)$ and $F_{n,\sinh }^{(j+1)} (x)$  (see Appendix~A). Then we group together the summands as follows
\begin{equation} \label{GrindEQ__29_}
\begin{split}
F_{n}^{(j+1)} (x)&=\sum _{p=0}^{M(j+1)-1} x^{p} \left(f_{n,\cos ,p}^{(j+1)} \cos \left(\frac{\pi n}{X} x\right)+f_{n,\sin ,p}^{(j+1)} \sin \left(\frac{\pi n}{X} x\right)\right)\\
&+\sum _{p=0}^{M(j)-1} x^{p} \left(f_{n,\cosh ,p}^{(j+1)} \cosh \left(\frac{\pi n}{X} x\right)+f_{n,\sinh ,p}^{(j+1)} \sinh \left(\frac{\pi n}{X} x\right)\right).
\end{split}
\end{equation}

To extract the coefficients of $x^p$ like
  $$f_{n,\cos ,p}^{(j+1)},\;f_{n,\sin ,p}^{(j+1)},\;p=0,1,...,M(j+1)-1$$
 and 
 $$f_{n,\cosh ,p}^{(j+1)},\;f_{n,\sinh ,p}^{(j+1)},\;p=0,1,...,M(j)-1$$
  in the polynomials $F_{n,\cos }^{(j+1)} (x)$, $F_{n,\sin }^{(j+1)} (x)$, $F_{n,\cosh }^{(j+1)} (x)$, $F_{n,\sinh }^{(j+1)} (x)$ (see Appendix A) noted above,  we can use  the function \verb|coeff|  with corresponding arguments in Maple. These coefficients are included in the main formulas of the proposed algorithm in Section~\ref{MakRom_section_6_1}. The expressions for these coefficients involve  only the algebraic operations and are represented through the corresponding quantities computed at previous steps of FD-method.

We require the polynomials at corresponding trigonometric functions and hyperbolic trigonometric functions to be equal on the both sides of equation  \eqref{GrindEQ__5_}, \eqref{GrindEQ__29_}. This requirement leads to the two recurrence systems \eqref{GrindEQ__39_}, \eqref{GrindEQ__40_} (with the initial conditions \eqref{GrindEQ__45_}, \eqref{GrindEQ__46_}) and \eqref{GrindEQ__47_}, \eqref{GrindEQ__48_} (with the initial conditions \eqref{GrindEQ__49_}, \eqref{GrindEQ__50_}) for the unknown coefficients of representation \eqref{GrindEQ__15_}.

The first system is the following:
\begin{equation} \label{GrindEQ__39_} 
\begin{split}
\left(\left(\left((t+4)b_{n,t+4}^{(j+1)} +4\frac{\pi n}{X} a_{n,t+3}^{(j+1)} \right)(t+3)-6\left(\frac{\pi n}{X} \right)^{2} b_{n,t+2}^{(j+1)} \right)(t+2)\right.&\\
\left.-4\left(\frac{\pi n}{X} \right)^{3} a_{n,t+1}^{(j+1)} \right)(t+1)=f_{n,\cos ,t}^{(j+1)} ,& 
\end{split}
\end{equation} 
\begin{equation} \label{GrindEQ__40_} 
\begin{split}
\left(\left(\left((t+4)a_{n,t+4}^{(j+1)} -4\frac{\pi n}{X} b_{n,t+3}^{(j+1)} \right)(t+3)-6\left(\frac{\pi n}{X} \right)^{2} a_{n,t+2}^{(j+1)} \right)(t+2)\right.&\\
\left.+4\left(\frac{\pi n}{X} \right)^{3} b_{n,t+1}^{(j+1)} \right)(t+1)=f_{n,\sin ,t}^{(j+1)} , & 
\end{split}
\end{equation} 
\[t=0,1,...,M\left(j+1\right)-4,\;\;\;\;\; j=0,1,2,...\] 
with the initial conditions
\begin{equation} \label{GrindEQ__45_} 
\begin{split}
a_{n,M(j+1)-s}^{(j+1)} &=\sum _{k=0}^{s}(-1)^{\left[\kern-0.15em\left[\left. \left. \frac{k}{2} \right]\kern-0.15em\right]\right. \right. +1} f_{n,\chi (k),M(j+1)-s-1+k}^{(j+1)}\\
&\times  \frac{(2k+1)(M(j+1)-1)^{k} }{2^{2+k} \left(M(j+1)-s+\left[\kern-0.35em\left[\left. \left. \frac{k+\delta _{2,s} }{2} \right]\kern-0.35em\right]\right. \right. \right)} \left(\frac{X}{\pi n} \right)^{3+k} , 
\end{split}
\end{equation} 
\begin{equation} \label{GrindEQ__46_} 
\begin{split}
b_{n,M(j+1)-s}^{(j+1)} &=\sum _{k=0}^{s}(-1)^{\left[\kern-0.15em\left[\left. \left. \frac{k+1}{2} \right]\kern-0.15em\right]\right. \right. } f_{n,\chi (k+1),M(j+1)-s-1+k}^{(j+1)}\\
&\times  \frac{(2k+1)(M(j+1)-1)^{k} }{2^{2+k} \left(M(j+1)-s+\left[\kern-0.35em\left[\left. \left. \frac{k+\delta _{2,s} }{2} \right]\kern-0.35em\right]\right. \right. \right)} \left(\frac{X}{\pi n} \right)^{3+k} , 
\end{split}
\end{equation} 
where $s=0,1,2$, $\chi (2p)=\cos$, $\chi (2p+1)=\sin$, $p=0,1$, $j=0,1,2,...$.
Here and below $\delta _{t,s} $ denotes the Kronecker delta, $\left[\kern-0.15em\left[y\right]\kern-0.15em\right]$ is the greatest integer less than or equal to a real number $y$ (in the computer algebra system Maple $\left[\kern-0.15em\left[y\right]\kern-0.15em\right]$  is the function \verb|floor(y)|). 

The second system is the following:
\begin{equation} \label{GrindEQ__47_} 
\begin{split}
\left(\left(\left((t+4)d_{n,t+4}^{(j+1)} +4\frac{\pi n}{X} c_{n,t+3}^{(j+1)} \right)(t+3)+6\left(\frac{\pi n}{X} \right)^{2} d_{n,t+2}^{(j+1)} \right)(t+2)\right.&\\
\left.+4\left(\frac{\pi n}{X} \right)^{3} c_{n,t+1}^{(j+1)} \right)(t+1)=f_{n,\cosh ,t}^{(j+1)} ,& 
\end{split}
\end{equation} 
\begin{equation} \label{GrindEQ__48_} 
\begin{split}
\left(\left(\left((t+4)c_{n,t+4}^{(j+1)} +4\frac{\pi n}{X} d_{n,t+3}^{(j+1)} \right)(t+3)+6\left(\frac{\pi n}{X} \right)^{2} c_{n,t+2}^{(j+1)} \right)(t+2)\right.&\\
\left.+4\left(\frac{\pi n}{X} \right)^{3} d_{n,t+1}^{(j+1)} \right)(t+1)=f_{n,\sinh ,t}^{(j+1)} ,&
\end{split} 
\end{equation}
\[t=0,1,...,M(j)-4,\;\;\;\;\; j=1,2,...\] 
with the initial conditions
\begin{equation} \label{GrindEQ__49_} 
\begin{split}
c_{n,M(j)-s}^{(j+1)} &=\sum _{k=0}^{s}(-1)^{\left[\kern-0.15em\left[\frac{k+1}{2} \right]\kern-0.15em\right]} f_{n,\chi h(k),M(j)-s-1+k}^{(j+1)}\\
& \times \frac{(2k+1)(M(j)-1)^{k} }{2^{2+k} \left(M(j)-s+\left[\kern-0.35em\left[\left. \left. \frac{k+\delta _{2,s} }{2} \right]\kern-0.35em\right]\right. \right. \right)} \left(\frac{X}{\pi n} \right)^{3+k}, 
\end{split}
\end{equation} 
\begin{equation} \label{GrindEQ__50_} 
\begin{split}
d_{n,M(j)-s}^{(j+1)} &=\sum _{k=0}^{s}(-1)^{\left[\kern-0.15em\left[\frac{k+1}{2} \right]\kern-0.15em\right]} f_{n,\chi h(k+1),M(j)-s-1+k}^{(j+1)}\\
&\times\frac{(2k+1)(M(j)-1)^{k} }{2^{2+k} \left(M(j)-s+\left[\kern-0.35em\left[\left. \left. \frac{k+\delta _{2,s} }{2} \right]\kern-0.35em\right]\right. \right. \right)} \left(\frac{X}{\pi n} \right)^{3+k} , 
\end{split}
\end{equation} 
where $s=0,1,2$, $\chi h (2p)=\cosh$, $\chi h (2p+1)=\sinh$, $p=0,1$, $j=1,2,...$.

 Let us introduce the column vectors:
\begin{equation} \label{GrindEQ__55_} 
\begin{split}
\vec{Z}_{n}^{[a,b]} (p,j)&=\left[a_{n,M(j+1)-p}^{(j+1)} ,b_{n,M(j+1)-p}^{(j+1)} \right]^{T} ,\\
\vec{F}_{n}^{[a,b]} (p+3,j)&=\frac{1}{4(M(j+1)-p-3)} \left(\frac{X}{\pi n} \right)^{3}\\
&\times \left[-f_{n,\cos ,M(j+1)-p-4}^{(j+1)} ,f_{n,\sin ,M(j+1)-p-4}^{(j+1)} \right]^{T} ,
\end{split}
\end{equation} 
\begin{equation} \label{GrindEQ__56_} 
\begin{split}
\vec{Z}_{n}^{[c,d]} (s,j)&=\left[c_{n,M(j)-s}^{(j+1)} ,d_{n,M(j)-s}^{(j+1)} \right]^{T} ,\\
\vec{F}_{n}^{[c,d]} (s+3,j)&=\frac{1}{4(M(j)-s-3)} \left(\frac{X}{\pi n} \right)^{3}\\
&\times \left[f_{n,\cosh ,M(j)-s-4}^{(j+1)} ,f_{n,\sinh ,M(j)-s-4}^{(j+1)} \right]^{T} ,
\end{split}
\end{equation} 
and denote the matrices:
\begin{equation} \label{GrindEQ__56_1_} 
\begin{split}
D_{n,1,1}^{[a,b]} (p,j)&=\frac{3X}{2\pi n} (M(j+1)-p-2)\left[\begin{array}{cc} {0} & {-1} \\ {1} & {0} \end{array}\right],\\
D_{n,1,2}^{[a,b]} (p,j)&=\left(\frac{X}{\pi n} \right)^{2} (M(j+1)-p-1)(M(j+1)-p-2)\left[\begin{array}{cc} {1} & {0} \\ {0} & {1} \end{array}\right],\\
D_{n,1,3}^{[a,b]} (p,j)&=-\frac{1}{4} \left(\frac{X}{\pi n} \right)^{3} (M(j+1)-p)(M(j+1)-p-1)\\
\end{split}
\end{equation} 
\begin{equation*} 
\begin{split}
&\times (M(j+1)-p-2)\left[\begin{array}{cc} {0} & {-1} \\ {1} & {0} \end{array}\right] 
\end{split}
\end{equation*} 
for $j=0,1,...$, and
\begin{equation} \label{GrindEQ__56_2_} 
\begin{split}
&D_{n,1,1}^{[c,d]} (s,j)=\left[\begin{array}{cc} {1} & {0} \\ {0} & {-1} \end{array}\right]D_{1,1}^{[a,b]} (s,j-1),\\
&D_{n,1,2}^{[c,d]} (s,j)=-D_{1,2}^{[a,b]} (s,j-1),\\
&D_{n,1,3}^{[c,d]} (s,j)=-\left[\begin{array}{cc} {1} & {0} \\ {0} & {-1} \end{array}\right]D_{1,3}^{[a,b]} (s,j-1)
\end{split}
\end{equation} 
for $j=1,2,...$. Here by $T$ we denote the transpose of a row vector.
For some fixed $j$ the systems of equations \eqref{GrindEQ__39_}--\eqref{GrindEQ__46_} and \eqref{GrindEQ__47_}--\eqref{GrindEQ__50_} can be expressed in vector form as linear inhomogeneous third-order difference equations with variable matrix coefficients
\begin{equation} \label{GrindEQ__57_} 
\begin{split}
\vec{Z}_{n}^{\nu} (p+3,j)&=D_{n,1,1}^{\nu} (p,j)\vec{Z}_{n}^{\nu} (p+2,j)+D_{n,1,2}^{\nu} (p,j)\vec{Z}_{n}^{\nu} (p+1,j)\\
&+D_{n,1,3}^{\nu} (p,j)\vec{Z}_{n}^{\nu} (p,j)+\vec{F}_{n}^{\nu} (p+3,j), 
\end{split}
\end{equation} 
where
$$p=0,1,...,M\left(j+1\right)-4,\;\;\;j=0,1,...\;\;\;\text{with}\;\;\;\nu=[a,b]$$
and 
$$p=0,1,...,M\left(j\right)-4,\;\;\;j=1,2,...\;\;\;\text{with}\;\;\;\nu=[c,d],$$
with the initial value vectors (see \eqref{GrindEQ__45_}, \eqref{GrindEQ__46_} and \eqref{GrindEQ__49_}, \eqref{GrindEQ__50_})
\begin{equation} \label{GrindEQ__58_} 
\vec{Z}_{n}^{\nu} (k,j),\; \; k=0,1,2,\; \; \;\nu=[a,b],\,[c,d].
\end{equation} 

Let us define the column vectors 
\begin{equation} \label{GrindEQ__61_} 
\begin{split}
\mathbf{Z}_{n}^{\nu} (p,j)&=\left[\vec{Z}_{n}^{\nu} (p+2,j),\vec{Z}_{n}^{\nu} (p+1,j),\vec{Z}_{n}^{\nu} (p,j)\right]^{T} ,\\
 \mathbf{F}_{n}^{\nu} (p)&=\left[\vec{F}_{n}^{\nu} (p+3,j),\vec{0},\vec{0}\right]^{T} ,
\end{split}
\end{equation} 
and the block matrix
\begin{equation} \label{GrindEQ__61_0_} 
\begin{split}
 \mathbf{D}_{n}^{\nu} (p,j)&=\left[\begin{array}{ccc} {D_{n,1,1}^{\nu} (p,j)} & {D_{n,1,2}^{\nu} (p,j)} & {D_{n,1,3}^{\nu} (p,j)} \\ {E} & {O} & {O} \\ {O} & {E} & {O} \end{array}\right]\\
 &=\left(D_{n,l_{1} ,l_{2} }^{\nu} (p,j)\right)_{l_{1} ,l_{2} =\overline{1,3}} 
\end{split}
\end{equation} 
with
\begin{equation} \label{GrindEQ__61_1_} 
\begin{split}
D_{n,2,1}^{\nu} (p,j)&=D_{n,3,2}^{\nu} (p,j)=E,\\
D_{n,2,2}^{\nu} (p,j)&=D_{n,2,3}^{\nu} (p,j)=D_{n,3,1}^{\nu} (p,j)=D_{n,3,3}^{\nu} (p,j)=O.
\end{split}
\end{equation} 
Here $E$ and $O$ are respectively the unit and zero $(2\times 2)$--matrices, $\vec{0}=\left[0,0\right]^{T} $ is the null column vector. 
We can rewrite \eqref{GrindEQ__57_}, \eqref{GrindEQ__58_} as
\begin{equation} \label{GrindEQ__62_} 
\begin{split}
&\mathbf{Z}_{n}^{\nu} (p+1,j)=\mathbf{D}_{n}^{\nu} (p,j)\mathbf{Z}_{n}^{\nu} (p,j)+\mathbf{F}_{n}^{\nu} (p,j),\\
&\mathbf{Z}_{n}^{\nu} (0,j)=\left[\vec{Z}_{n}^{\nu} (2,j),\vec{Z}_{n}^{\nu} (1,j),\vec{Z}_{n}^{\nu} (0,j)\right]^{T} ,
\end{split}
\end{equation} 
where 
$$p=0,1,...,M\left(j+1\right)-4,\;\;\;j=0,1,...\;\;\;\text{with}\;\;\;\nu=[a,b]$$
and 
$$p=0,1,...,M\left(j\right)-4,\;\;\;j=1,2,...\;\;\;\text{with}\;\;\;\nu=[c,d].$$

The solution of \eqref{GrindEQ__62_} can then be given in matrix form by
\begin{equation} \label{GrindEQ__63_} 
\begin{split}
&\mathbf{Z}_{n}^{\nu} (p+1,j)=\mathbf{D}_{n}^{\nu} (p,j)\mathbf{D}_{n}^{\nu} (p-1,j)\cdot \ldots \cdot \mathbf{D}_{n}^{\nu} (0,j)\mathbf{Z}_{n}^{\nu} (0,j)+\\
&+\sum _{s=1}^{p}\mathbf{D}_{n}^{\nu} (p,j)\mathbf{D}_{n}^{\nu} (p-1,j)\cdot \ldots \cdot \mathbf{D}_{n}^{\nu} (p-s+1,j)\mathbf{F}_{n}^{\nu} (p-s,j) +\mathbf{F}_{n}^{\nu} (p,j).
\end{split}
\end{equation} 
Hence the solution of \eqref{GrindEQ__57_} is  
\begin{equation} \label{GrindEQ__64_} 
\begin{split}
&\vec{Z}_{n}^{\nu} (p+3,j)=\sum_{\begin{subarray}{l} {1\le l_{p} ,l_{p-1} ,...,l_{0} \le 3,} \\ {l_{p+1} =1} \end{subarray}}
\prod _{s=0}^{p}D_{n,l_{s+1} ,l_{s} }^{\nu} (s,j)  \vec{Z}_{n}^{\nu} (3-l_{0} ,j)\\
&+\sum _{s=1}^{p} \sum _{\begin{subarray}{l} {1\le l_{s-1} ,...,l_{1} \le 3,} \\ {l_{0} =1,\; l_{s} =1} \end{subarray}}\prod _{k=1}^{s}D_{n,l_{k} ,l_{k-1} }^{\nu} (p-s+k,j)  \vec{F}_{n}^{\nu} (p-s+3,j)\\
&+\vec{F}_{n}^{\nu} (p+3,j),
\end{split}
\end{equation}
where 
$$p=0,1,...,M\left(j+1\right)-4,\;\;\;j=0,1,...\;\;\;\text{with}\;\;\;\nu=[a,b]$$
and 
$$p=0,1,...,M\left(j\right)-4,\;\;\;j=1,2,...\;\;\;\text{with}\;\;\;\nu=[c,d].$$

\begin{remark}\label{remark_MakRom_1_0} Note, if the upper bound of summation index is less than the lower bound of summation index, then the sum is the empty sum, with the value $0$.
In \eqref{GrindEQ__64_} and below in \eqref{GrindEQ__65_1_}, \eqref{GrindEQ__65_2_} we multiply the matrices $D_{n,l_{s+1} ,l_{s} }^{\nu} (s,j)$ from the right to the left, i.e.,
\begin{equation*} 
\begin{split}
&\prod _{s=0}^{p}D_{n,l_{s+1} ,l_{s} }^{\nu} (s,j) \vec{Z}_{n}^{\nu} (3-l_{0} ,j)\\
&=D_{n,1,l_{p} }^{\nu} (p,j)\cdot \left(D_{n,l_{p} ,l_{p-1} }^{\nu} (p-1,j)\cdot \right. \left(\ldots  \right. \\
& \ldots \left. \left. \cdot\left(D_{n,l_{2} ,l_{1} }^{\nu} (1,j)\left(D_{n,l_{1} ,l_{0} }^{\nu} (0,j)\vec{Z}_{n}^{\nu} (3-l_{0} ,j)\right)\right)\cdots \right)\right).
\end{split}
\end{equation*} 
In \eqref{GrindEQ__64_} and below in \eqref{GrindEQ__65_1_}, \eqref{GrindEQ__65_2_} the notation
$$\sum _{\begin{subarray}{l} {\; \; \; 1\le l_{p} ,l_{p-1} ,...,l_{0} \le 3,} \\ {\, \, \, \, \, \, l_{p+1} =1} \end{subarray}}$$
denotes the following expression
$$\sum _{l_{p+1} =1}^{1}\sum _{l_{p} =1}^{3}\sum _{l_{p-1} =1}^{3}\cdots    \sum _{l_{2} =1}^{3}\sum _{l_{1} =1}^{3}\sum _{l_{0} =1}^{3}.$$
\end{remark}

Returning to the substitution \eqref{GrindEQ__55_}, \eqref{GrindEQ__56_}, from \eqref{GrindEQ__64_} we
obtain the following column vectors: 
\begin{equation} \label{GrindEQ__65_1_} 
\begin{split}
&\left[\!\!\begin{array}{l} {a_{n,M(j+1)-p-3}^{(j+1)} } \\ {b_{n,M(j+1)-p-3}^{(j+1)} } \end{array}\!\!\right]\!\!=\!\sum _{\begin{subarray}{l} {1\le l_{p} ,l_{p-1} ,...,l_{0} \le 3,} \\ {l_{p+1} =1} \end{subarray}}\prod _{s=0}^{p}\!D_{n,l_{s+1} ,l_{s} }^{[a,b]} (s,j)  \!\!\left[\!\!\begin{array}{l} {a_{n,M(j+1)-3+l_{0} }^{(j+1)} } \\ {b_{n,M(j+1)-3+l_{0} }^{(j+1)} } \end{array}\!\!\right]\\
&+\frac{1}{4} \left(\frac{X}{\pi n} \right)^{3} \sum _{s=1}^{p}\frac{1}{M(j+1)+s-p-3}  \\
&\times\sum _{\begin{subarray}{l} {1\le l_{s-1} ,...,l_{1} \le 3,} \\ {l_{0} =1,\, \, l_{s} =1} \end{subarray}}\prod _{k=1}^{s}\!D_{n,l_{k} ,l_{k-1} }^{[a,b]} (p-s+k,j) \! \left[\!\begin{array}{l} {-f_{n,\cos ,M(j+1)-p+s-4}^{(j+1)} } \\ {f_{n,\sin ,M(j+1)-p+s-4}^{(j+1)} } \end{array}\!\right]\\
& +\frac{1}{4(M(j+1)-p-3)} \left(\frac{X}{\pi n} \right)^{3}  \left[\begin{array}{l} {-f_{n,\cos ,M(j+1)-p-4}^{(j+1)} } \\ {f_{n,\sin ,M(j+1)-p-4}^{(j+1)} } \end{array}\right],
\end{split}
\end{equation} 
$$p=0,1,...,M\left(j+1\right)-4,\;\;\;j=0,1,2,...$$
and
\begin{equation} \label{GrindEQ__65_2_} 
\begin{split}
&\left[\begin{array}{l} {c_{n,M(j)-p-3}^{(j+1)} } \\ {d_{n,M(j)-p-3}^{(j+1)} } \end{array}\right] =\sum _{\begin{subarray}{l} {1\le l_{p} ,l_{p-1} ,...,l_{0} \le 3,} \\ {l_{p+1} =1} \end{subarray}}\prod _{s=0}^{p}D_{n,l_{s+1} ,l_{s} }^{[c,d]} (s,j)  \left[\begin{array}{l} {c_{n,M(j)-3+l_{0} }^{(j+1)} } \\ {d_{n,M(j)-3+l_{0} }^{(j+1)} } \end{array}\right]\\
&+\frac{1}{4} \left(\frac{X}{\pi n} \right)^{3} \sum _{s=1}^{p}\frac{1}{4(M(j)+s-p-3)}\\
&\times \sum _{\begin{subarray}{l} {1\le l_{s-1} ,...,l_{1} \le 3,} \\ {l_{0} =1,\, \, l_{s} =1} \end{subarray}}\prod _{k=1}^{s}D_{n,l_{k} ,l_{k-1} }^{[c,d]} (p-s+k,j)  \; \left[\begin{array}{l} {f_{n,\cosh ,M(j)-p+s-4}^{(j+1)} } \\ {f_{n,\sinh ,M(j)-p+s-4}^{(j+1)} } \end{array}\right]\\
&+\frac{1}{4(M(j)-p-3)} \left(\frac{X}{\pi n} \right)^{3} \left[\begin{array}{l} {f_{n,\cosh ,M(j)-p-4}^{(j+1)} } \\ {f_{n,\sinh ,M(j)-p-4}^{(j+1)} } \end{array}\right] ,
\end{split}
\end{equation} 
$$p=0,1,...,M\left(j\right)-4,\;\;\;j=1,2,....$$
Namely, we obtain the explicit recursive formulas for the coefficients in \eqref{GrindEQ__15_} 
\[a_{n,M(j+1)-p-3}^{(j+1)} ,\;\;\;\;b_{n,M(j+1)-p-3}^{(j+1)} ,\;\;\;\; p=0,1,...,M\left(j+1\right)-4,\;\;\;\; j=0,1,2,...\] 
and
\[c_{n,M(j)-p-3}^{(j+1)},\;\;\;\;d_{n,M(j)-p-3}^{(j+1)} ,\;\;\;\; p=0,1,...,M\left(j\right)-4,\;\;\;\; j=1,2,...\]
which are the corresponding elements of the vectors \eqref{GrindEQ__65_1_} and \eqref{GrindEQ__65_2_}.
These coefficients are represented recursively through the corresponding coefficients and quantities computed at previous steps of FD-method as well as through the coefficients of the polynomials \eqref{GrindEQ__3_}.

Substituting the representation \eqref{GrindEQ__15_} in the boundary conditions \eqref{GrindEQ__6_}, 
we obtain the nonhomogeneous system of linear algebraic equations for the coefficients
$b_{n,0}^{(j+1)}$, $c_{n,0}^{(j+1)}$, $d_{n,0}^{(j+1)}$, $j=0,1,2,...$. The solution of this system is:
\begin{equation} \label{GrindEQ__66_1_} 
\begin{split}
b_{n,0}^{(j+1)} &=\frac{X}{\pi n} \left(a_{n,1}^{(j+1)} +c_{n,1}^{(j+1)} +\frac{X}{\pi n} \left(b_{n,2}^{(j+1)} +d_{n,2}^{(j+1)} \right)\right),\; c_{n,1}^{(1)} =d_{n,2}^{(1)} =0,\\
d_{n,0}^{(j+1)} &=-b_{n,0}^{(j+1)} ,\\
c_{n,0}^{(j+1)}\! &=\!-\frac{1}{\sinh \left(\pi n\right)} \left(\!\sum _{t=0}^{M(j+1)}X^{t} b_{n,t}^{(j+1)} \! \cos \left(\pi n\right)+ \sum _{t=0}^{M(j)}X^{t} d_{n,t}^{(j+1)} \! \cosh \left(\pi n\right)\!\right)\\
&-\sum _{t=1}^{M(j)}X^{t} c_{n,t}^{(j+1)}  ,\;\;\;j=0,1,2,...
\end{split}
\end{equation}
The constant $a_{n,0}^{(j+1)}$ is calculated by the formula 
\begin{equation} \label{GrindEQ__66_2_} 
\begin{split}
a_{n,0}^{(j+1)}&=-\frac{2}{X} \left(\sum _{t=1}^{M(j+1)}\left(\beta _{n,t} b_{n,t}^{(j+1)} +\alpha _{n,t} a_{n,t}^{(j+1)} \right) \right. \\
&+\left. \sum _{t=0}^{M(j)}\left(\eta _{n,t} d_{n,t}^{(j+1)} +\mu _{n,t} c_{n,t}^{(j+1)} \right) \right),\;\;\;\;
j=0,1,2,...
\end{split}
\end{equation}
obtained from the orthogonality condition \eqref{GrindEQ__8_} and the formulas \eqref{GrindEQ__66_1_}.
Here the notations $\alpha _{n,t}$, $\beta _{n,t}$, $\eta _{n,t}$, $\mu _{n,t}$ are used which are exactly calculated in the case $t=0,1,...$. The analytical expressions for these notations are given in Appendix~B (see also \cite{GradshteynRyzhik2014}).

Using Lemma~\ref{MakRom_lemma_2}, from \eqref{GrindEQ__7_} we obtain the formula for the corrections of the eigenvalues $\lambda _{n}^{(j+1)}$ which is given in Appendix~C. 

\begin{remark}\label{remark_MakRom_1}
Analytical expressions for the approximations $\mathop{\lambda_{n} }\limits^{m}$, $\mathop{u_{n}}\limits^{m}  (x)$ of rank $m$ (according to the FD-method) to the exact eigenpairs $\lambda _{n}$, $u_{n}(x)$ are the results of the execution of a numerical algorithm from Section~\ref{MakRom_section_6_1}.
These expressions analytically depend on the eigenpair index number $n$ and on the input data of the problem under consideration \eqref{GrindEQ__1_}--\eqref{GrindEQ__3_} as $X$, $r_{0}$, $r_{1}$, $r_{2}$, $A_{l}$ $(l=0,1,...,r_{0})$, $B_{l}$ $(l=0,1,...,r_{1})$, $C_{l}$ $(l=0,1,...,r_{2} )$.
In order to find numerically result for a given value of $n$, we calculate the corresponding approximations $\mathop{\lambda_{n} }\limits^{m}$, $\mathop{u_{n}}\limits^{m}  (x)$, by substituting in the obtained analytical expressions the value $n$ and the numerical values of the parameters entering into the input data, if there are any.
\end{remark}

\newpage

\section{Numerical algorithm}\label{MakRom_section_6_1}%6_1

\textbf{Data:} \textbf{choose} a set of input values: $X,$ $r_{0} ,r_{1} ,r_{2} ,$ $A_{l} (l=0,1,...,r_{0} ),$ \newline
$B_{l} (l=0,1,...,r_{1} ),$ $C_{l} (l=0,1,...,r_{2} )$, $m$ $(m \geq 1)$ (see problem \eqref{GrindEQ__1_}--\eqref{GrindEQ__3_});

\textbf{Result:} $\mathop{\lambda_{n}}\limits^{m}$, $\mathop{u_{n}}\limits^{m}(x)$;

\begin{enumerate}
\item[1)]  \textbf{set input parameters} for a recurrence procedure
$r=\max \left\{r_{0} ,r_{1} ,r_{2} \right\}$, $a_{n,0}^{(0)}=\sqrt{\frac{2}{X} }$, $b_{n,0}^{(0)} =d_{n,0}^{(0)} =c_{n,0}^{(0)} =0$, $\lambda _{n}^{(0)} =\frac{\left(n\pi \right)^{4} }{X^{4} }$, $u_{n}^{(0)} (x)=a_{n,0}^{(0)} \sin \left(\frac{n\pi }{X} x\right)$;

\item[2)]  \textbf{compute} the values of ancillary quantities $\beta _{n,t}$, $\alpha _{n,t}$, $\eta _{n,s}$, $\mu _{n,s}$, 
 $(t=0,1,...,M(m)-1$, $s=0,1,...,M(m-1)-1)$ using Appendix~B;
 
\item[3)]  \textbf{compute} the correction $\lambda _{n}^{(1)}$ given in Appendix~C with $j=0$; 

\item[4)]  \textbf{compute} the functions $F_{n,\cos }^{(1)} (x)$, $F_{n,\sin }^{(1)} (x)$ given in Appendix~A with $j=0$ and \textbf{compute} the corresponding coefficients $f_{n,\cos ,p}^{(1)}$, $f_{n,\sin ,p}^{(1)}$ $(p=0,1,...,r)$ (see \eqref{GrindEQ__29_}); 

\item[5)]  \textbf{compute} the coefficients $a_{n,r+1-s}^{(1)}$, $b_{n,r+1-s}^{(1)}$ $(s=0,1,2)$ using \eqref{GrindEQ__45_}, \eqref{GrindEQ__46_} with $j=0$;

\item[6)] using \eqref{GrindEQ__56_1_}, \eqref{GrindEQ__61_1_} with $j=0$ \textbf{compute} the matrices 
$$D_{n,l_{s+1} ,l_{s} }^{[a,b]} (s,0)\;\; (l_{s} =1,2,3,\;s=0,1,...,r-3);$$
$$D_{n,l_{k} ,l_{k-1} }^{[a,b]} (p-s+k,0)\;\;(l_{k} =1,2,3,\;k=1,2,...,s,\; s=1,2,...,p,$$
$$p=0,1,...,r-3);$$ 
 
\item[7)]  \textbf{compute} the coefficients $a_{n,r-p-2}^{(1)}$, $b_{n,r-p-2}^{(1)}$ $(p=0,1,...,r-3)$ using \eqref{GrindEQ__65_1_} with $j=0$;

\item[8)]  \textbf{compute} the coefficients $b_{n,0}^{(1)}$, $c_{n,0}^{(1)}$, $d_{n,0}^{(1)}$, $a_{n,0}^{(1)} $ using \eqref{GrindEQ__66_1_}, \eqref{GrindEQ__66_2_} with $j=0$;

\item[9)]  \textbf{compute} the correction $u_{n}^{(1)} (x)$ using \eqref{GrindEQ__15_} with $j=0$;\newline
\textbf{if} $m>1$ \textbf{then}

\textbf{for $j$ from $1$ (with unit step) to $m-1$ do}

\item[10)] \textbf{compute} $\lambda _{n}^{(j+1)}$ given in Appendix~C;

\item[11)]  \textbf{compute} the functions $F_{n,\cos }^{(j+1)} (x)$, $F_{n,\sin }^{(j+1)} (x)$, $F_{n,\cosh }^{(j+1)} (x)$, $F_{n,\sinh }^{(j+1)} (x)$ given in Appendix~A and \textbf{compute} the corresponding coefficients $f_{n,\cos ,p}^{(j+1)}$, $f_{n,\sin ,p}^{(j+1)}$ $(p=0,1,...,M(j+1)-1)$ and $f_{n,\cosh ,p}^{(j+1)}$, $f_{n,\sinh ,p}^{(j+1)}$ $(p=0,1,...,M(j)-1)$ (see \eqref{GrindEQ__29_}); 

\item[12)]  \textbf{compute} the coefficients $a_{n,M(j+1)-s}^{(j+1)}$, $b_{n,M(j+1)-s}^{(j+1)}$, $c_{n,M(j)-s}^{(j+1)}$, $d_{n,M(j)-s}^{(j+1)}$ $(s=0,1,2)$ using \eqref{GrindEQ__45_}, \eqref{GrindEQ__46_}, \eqref{GrindEQ__49_}, \eqref{GrindEQ__50_};

\item[13)] using \eqref{GrindEQ__56_1_}, \eqref{GrindEQ__61_1_} \textbf{compute} the matrices 
$$D_{n,l_{s+1} ,l_{s} }^{[a,b]} (s,j)\;\; (l_{s} =1,2,3,\;s=0,1,...,M(j+1)-4),$$  
$$D_{n,l_{k} ,l_{k-1} }^{[a,b]} (p-s+k,j)\;\;(l_{k} =1,2,3,\;k=1,2,...,s,\; s=1,2,...,p,$$
$$p=0,1,...,M(j+1)-4),$$ 
$$D_{n,l_{s+1} ,l_{s} }^{[c,d]} (s,j)\;\;(l_{s} =1,2,3,\;s=0,1,...,M(j)-4),$$  
$$D_{n,l_{k} ,l_{k-1} }^{[c,d]} (p-s+k,j)\;\;(l_{k} =1,2,3,\;k=1,2,...,s,\; s=1,2,...,p,$$ $$p=0,1,...,M(j)-4);$$ 

\item[14)] \textbf{compute} the coefficients $a_{n,M(j+1)-p-3}^{(j+1)}$, $b_{n,M(j+1)-p-3}^{(j+1)}$ $(p=0,1,...,M(j+1)-4)$ using \eqref{GrindEQ__65_1_} and \textbf{compute} the coefficients $c_{n,M(j)-p-3}^{(j+1)}$, $d_{n,M(j)-p-3}^{(j+1)}$ $(p=0,1,...,M(j)-4)$ using \eqref{GrindEQ__65_2_};

\item[15)]  \textbf{compute} the coefficients $b_{n,0}^{(j+1)}$, $c_{n,0}^{(j+1)}$, $d_{n,0}^{(j+1)}$, $a_{n,0}^{(j+1)}$ using \eqref{GrindEQ__66_1_}, \eqref{GrindEQ__66_2_};

\item[16)]  \textbf{compute} the correction $u_{n}^{(j+1)} (x)$ using \eqref{GrindEQ__15_};

\textbf{end for;}

\textbf{end if;}

\item[17)]  \textbf{compute} the approximations $\mathop{u_{n}}\limits^{m}(x)$, $\mathop{\lambda_{n} }\limits^{m}$ using \eqref{GrindEQ__4_}.

\end{enumerate}

\newpage

\section{Numerical examples}\label{MakRom_section_7}%7

\begin{example}\label{example_MakRom_1}
We consider the eigenvalue problem \eqref{GrindEQ__1_}, \eqref{GrindEQ__2_} with $X=5$, and with potential coefficients in \eqref{GrindEQ__3_} equal to
\begin{equation} \label{GrindEQ__67_} 
q_{0} (x)=0.0001x^{4} -0.02,\;\;\;q_{1} (x)=-0.04x,\;\;\;q_{2} (x)=-0.02x^{2}.
\end{equation}

The computations of the exact eigenvalues $\lambda_n$ and eigenfunctions $u_n (x)$, and of their approximations $ \mathop{\lambda_{n}}\limits^{m}$ and $\mathop{u_{n}}\limits^{m}(x)$ obtained by FD-method of rank $m$,  have been done with the help of the computer algebra system Maple (\verb|Digits=300|).

The exact solution of the problem \eqref{GrindEQ__1_}, \eqref{GrindEQ__2_}, \eqref{GrindEQ__67_} is expressed in terms of confluent hypergeometric Kummer's functions $M\left(a,b,z\right)$ and $U\left(a,b,z\right)$ (see \cite[Chapter~13]{NIST2010})
\begin{equation} \label{GrindEQ__68_} 
\begin{split}
u_n &(x)\!=\!\frac{xe^{-{x^{2} \mathord{\left/ {\vphantom {x^{2}  20}} \right. \kern-\nulldelimiterspace} 20} } }{10^{{3\mathord{\left/ {\vphantom {3 4}} \right. \kern-\nulldelimiterspace} 4} } }\! \left(\!C_{1} M\!\left(\frac{3}{4} +\frac{5}{2} \sqrt{\lambda_n } ,\frac{3}{2} ,\frac{x^{2} }{10} \right)\right. \!+\!C_{2} M\!\left(\frac{3}{4}\!-\!\frac{5}{2} \sqrt{\lambda_n } ,\frac{3}{2} ,\frac{x^{2} }{10} \right)\\
&+\!C_{3} U\left(\frac{3}{4}\!+\!\frac{5}{2} \sqrt{\lambda_n } ,\frac{3}{2} ,\frac{x^{2} }{10} \right)\!+\!C_{4}\! \left. U\left(\frac{3}{4}\!-\!\frac{5}{2} \sqrt{\lambda_n } ,\frac{3}{2} ,\frac{x^{2} }{10} \right)\!\right).
\end{split}
\end{equation}
The constants $C_{k}$, $k=1,2,3,4$ and the eigenvalues $\lambda_n$ of the problem \eqref{GrindEQ__1_}, \eqref{GrindEQ__2_}, \eqref{GrindEQ__67_} can be found from the system of equations obtained by substituting \eqref{GrindEQ__68_} into the boundary conditions \eqref{GrindEQ__2_} with $X=5$.
Setting the determinant of this system equal to zero, we obtain a transcendental equation with respect to $\lambda_n$
\begin{equation} \label{GrindEQ__69_} 
\sqrt{\lambda_n } M\left(\frac{3}{4} +\frac{5}{2} \lambda_n ,\frac{3}{2} ,\frac{5}{2} \right)M\left(\frac{3}{4} -\frac{5}{2} \lambda_n ,\frac{3}{2} ,\frac{5}{2} \right)=0,\;\;\;\;n=1,2,...\; . 
\end{equation} 
Using the command \verb|fsolve| in Maple, from \eqref{GrindEQ__69_} we find the first eight smallest exact eigenvalues $\lambda _{n}$, $n=1,2,...,8$ of the problem under consideration, which are given in Table~\ref{table_MakRom_1}.

\begin{table}\label{table_MakRom_1}
\begin{center}
{\footnotesize{
\tabcolsep=0.11cm
\begin{tabular}{p{0.2in}  p{3.0in}} \hline 
$n$ & $\lambda _{n}$ \\ \hline 
1 & 0.2150508643697154969799099152379067104370468531017 \\ 
2 & 2.7548099346830341769807978567162582697401440677181 \\  
3 & 13.215351540558178725583135747686922412859870693289 \\ 
4 & 40.950819759161479687386430406406575369564981705685 \\  
5 & 99.053478063489519905364004392278295950913079764143 \\  
6 & 204.35573226825688655101487325028333585898058873698 \\ 
7 & 377.43042068923559313999045804039609985719687266348 \\ 
8 & 642.59086816966269512711767453062929784652398288740 \\ \hline 
\end{tabular}
\caption{Exact eigenvalues $\lambda _{n}$ for $n=1,2,...,8$ of the problem \eqref{GrindEQ__1_}, \eqref{GrindEQ__2_}, \eqref{GrindEQ__67_} from Example~\ref{example_MakRom_1}.}
}}
\end{center}
\end{table}

Moreover, the eigenvalues $\lambda _{n}$ of the given problem \eqref{GrindEQ__1_}, \eqref{GrindEQ__2_}, \eqref{GrindEQ__67_} are precisely the squares of the corresponding eigenvalues $L_n$ of the second-order Sturm--Liouville problem with Dirichlet boundary conditions (see \cite[Section~6]{69_MarlettaGreenberg1997}):
\begin{equation} \label{GrindEQ__70_} 
y''(x)+\left(L_n-0.01x^{2} \right)y(x)=0,\;x\in \; {\kern 1pt} (0,5),\; y(0)=y(5)=0. 
\end{equation} 
The exact solution of \eqref{GrindEQ__70_} is 
\begin{equation} \label{GrindEQ__71_} 
y(x)=\frac{xe^{-{x^{2} \mathord{\left/ {\vphantom {x^{2}  20}} \right. \kern-\nulldelimiterspace} 20} } }{10^{{3\mathord{\left/ {\vphantom {3 4}} \right. \kern-\nulldelimiterspace} 4} } }\left(C_{1} M\left(\frac{3}{4} -\frac{5}{2} L_n,\frac{3}{2} ,\frac{1}{10} x^{2} \right)+C_{2} U\left(\frac{3}{4} -\frac{5}{2} L_n,\frac{3}{2} ,\frac{1}{10} x^{2} \right)\right).
\end{equation}
The constants $C_{1} $, $C_{2} $  and the eigenvalues $L_n$ of the problem \eqref{GrindEQ__70_} can be found from the system of equations obtained by substituting \eqref{GrindEQ__71_} into the boundary conditions $y(0)=y(5)=0$.
Setting the determinant of this system equal to zero, we obtain a transcendental equation with respect to $L_n$
\begin{equation} \label{GrindEQ__72_} 
M\left(\frac{3}{4} -\frac{5}{2} L_n,\frac{3}{2} ,\frac{5}{2} \right)=0,\;\;\;\;n=1,2,...\;. 
\end{equation} 
Exact eigenvalues $L_n$ can be found from \eqref{GrindEQ__72_} using the command \verb|fsolve| in Maple. The squares of the eigenvalues $L_n$, i.e. the values $\lambda _{n}=\left(L_n\right)^{2}$, are given in Table~\ref{table_MakRom_1}. 

Proposed in this paper symbolic algorithm of the FD-method was applied to compute the approximate solution of the problem \eqref{GrindEQ__1_}, \eqref{GrindEQ__2_}, \eqref{GrindEQ__67_}. The approximate eigenpairs $\mathop{u _{n}}\limits^{m} (x)$, $\mathop{\lambda_{n} }\limits^{m}$ were computed exactly as analytical expressions with respect to the index number $n$ and the input data of the problem under consideration, i.e., we had no rounding errors (see Remark~\ref{remark_MakRom_1}).
Below we give the eigenvalue corrections $\lambda _{n}^{(j+1)}$ for some first iteration steps $j$ of FD-method:
\begin{equation*} \label{GrindEQ__73_} 
\begin{split}
\lambda _{n}^{(0)} &=\frac{\left(n\pi \right)^{4} }{625} ,\lambda _{n}^{(1)} =\frac{\left(n\pi \right)^{2} }{150} +\frac{1}{400} -\frac{1}{16\left(n\pi \right)^{2} } +\frac{3}{32\left(n\pi \right)^{4} } ,\\
\lambda _{n}^{(2)} &=-\frac{1}{360} +\frac{1}{224\left(n\pi \right)^{2} } -\frac{173}{384\left(n\pi \right)^{4} } +\frac{9075}{3584\left(n\pi \right)^{6} } -\frac{625\coth \left(n\pi \right)}{64\left(n\pi \right)^{7} } \\
&+\frac{28775}{512\left(n\pi \right)^{8} } -\frac{556875}{2048\left(n\pi \right)^{10} } +\frac{804375}{2048\left(n\pi \right)^{12} } ,
\end{split}
\end{equation*}
\begin{equation*}
\begin{split}
\lambda _{n}^{(3)} &=-\frac{25}{4472832\left(e^{2\pi n} -1\right)^{2} } \left(\left(\frac{11748744}{\left(\pi n\right)^{8} } -\frac{338782925}{\left(\pi n\right)^{12} } \right)\left(e^{4\pi n} +1\right)\right.\\
\end{split}
\end{equation*}
\begin{equation*}
\begin{split}
&\left.+\left(-\frac{38057488}{\left(\pi n\right)^{8} } +\frac{814065850}{\left(\pi n\right)^{12} } \right)e^{2\pi n} \right)+\frac{625}{3072} \left(\frac{96}{\left(\pi n\right)^{7} } +\frac{400}{\left(\pi n\right)^{9} } \right.
\end{split}
\end{equation*}
\begin{equation*}
\begin{split}
&\left.-\frac{14250}{\left(\pi n\right)^{11} } -\frac{5625}{\left(\pi n\right)^{13} } +\frac{{\rm 1203750}}{\left(\pi n\right)^{15} } \right)\cosh \left(\pi n\right)-\frac{5}{672\left(\pi n\right)^{2} } +\frac{155}{2688\left(\pi n\right)^{4} }\\
&-\frac{51325}{39424\left(\pi n\right)^{6} } +\frac{{\rm 232994375}}{{\rm 315392}\left(\pi n\right)^{10} } -\frac{{\rm 24718046875}}{{\rm 229376}\left(\pi n\right)^{14} } +\frac{{\rm 54692015625}}{{\rm 65536}\left(\pi n\right)^{16} }\\
&-\frac{{\rm 843444140625}}{{\rm 131072}\left(\pi n\right)^{18} } +\frac{{\rm 675516796875}}{{\rm 65536}\left(\pi n\right)^{20} }.
\end{split}
\end{equation*}
Here the command \verb|combine(,trig)| in Maple was used to rewrite the eigenvalue corrections $\lambda _{n}^{(j+1)}$  in a compact form.
We give the coefficients of \eqref{GrindEQ__15_} only for first two steps of the FD-method with $j=-1,0$ because they are too large
\begin{equation*}
\begin{split}
&a_{n,0}^{(0)}=\sqrt{\frac{2}{5} },\;\;\;b_{n,0}^{(0)} =c_{n,0}^{(0)} =d_{n,0}^{(0)} =0\; \text{(step 1)};\\
&b_{n,5}^{(1)} =-\frac{\sqrt{10} }{8000\left(\pi n\right)^{3} } ,\;\;\;a_{n,4}^{(1)} =\frac{3\sqrt{10} }{640\left(\pi n\right)^{4} } ,\;\;\;
b_{n,3}^{(1)} =\frac{\sqrt{10} }{8\pi n} \left(-\frac{1}{75} +\frac{5}{8\left(\pi n\right)^{4} } \right),\\
&a_{n,5}^{(1)} =b_{n,4}^{(1)} =a_{n,3}^{(1)} =0\;\text{(step 5)};\\
&a_{n,2}^{(1)} =\frac{\sqrt{10} }{16\left(\pi n\right)^{2} } \left(\frac{1}{5} -\frac{75}{8\left(\pi n\right)^{4} } \right),\;\;\;b_{n,1}^{(1)} =\frac{\sqrt{10} }{8\pi n} \left(\frac{1}{3} +\frac{5}{8\left(\pi n\right)^{2} } -\frac{25}{8\left(\pi n\right)^{4} } \right),\\
&b_{n,2}^{(1)} =a_{n,1}^{(1)} =0\;\text{(step 7)};\\
&c_{n,0}^{(1)} =-\frac{125\sqrt{10} \cos \left(\pi n\right)}{16\left(\pi n\right)^{5} \sinh \left(\pi n\right)} ,\;a_{n,0}^{(1)} \!=\!\frac{5\sqrt{10} }{128\left(\pi n\right)^{2} } \left(\!-\frac{8}{3}\! -\!\frac{7}{\left(\pi n\right)^{2} } \!+\!\frac{125}{\left(\pi n\right)^{4} }\! -\!\frac{525}{\left(\pi n\right)^{6} }\! \right),\\
&b_{n,0}^{(1)} =d_{n,0}^{(1)} =0 \;\text{(step 8)}.
\end{split}
\end{equation*}
Here in parentheses the step of the numerical algorithm from Section~\ref{MakRom_section_6_1} is given on which these coefficients are calculated.

Using the proposed symbolic algorithm of the FD-method of rank $m=0, 1, 2, ..., 20$ the approximations $\mathop{\lambda_{n} }\limits^{m}$ to the first eight eigenvalues $\lambda_n$ with $n=1,2,...,8$ were calculated by substitution the value of $n$ into the analytical expressions for $\mathop{\lambda_{n} }\limits^{m}$.
Figure~1 shows the broken logarithmic line graphs which were created connecting the data points $\left(m;\ln \left(\Delta _{n}^{FD} (m)\right)\right)$ by lines, where $\Delta _{n}^{FD} (m)= |\lambda _{n} -\mathop{\lambda _{n} }\limits^{m} |$ are the absolute errors of the approximations $\mathop{\lambda_{n} }\limits^{m}$ of rank $m$ to the exact eigenvalues $\lambda_{n}$ with the number $n$. Figure~1 illustrates the exponential convergence of the proposed approach for the problem \eqref{GrindEQ__1_}, \eqref{GrindEQ__2_}, \eqref{GrindEQ__67_}. The sufficient convergence condition \eqref{GrindEQ__26_} (with $\omega =0.2$) is fulfilled for $n\ge 3$, but, as we can see in Figure~1, the method converges for $n=1,2$ too, i.e., the conditions of Theorem~\ref{MakRom_theorem_3} are rough and can be improved.

\begin{figure}\label{MakRom_fig_1}
\begin{center}
\includegraphics[width=80mm]{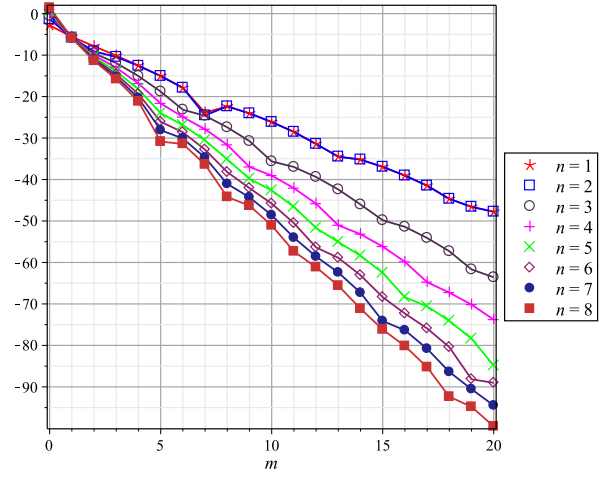}
\caption{Graphs for broken lines joining the data points $\left(m;\ln \left(\Delta _{n}^{FD} (m)\right)\right)$ for Example~\ref{example_MakRom_1}. These graphs illustrates the behaviour of the absolute errors $\Delta _{n}^{FD} (m)=|\lambda _{n} -\mathop{\lambda _{n} }\limits^{m} |$ of FD-method for $m=0,1,2,...,20$, i.e., they illustrate the exponential convergence of the FD-method with respect to the rank of FD-method $m$ for the eigenpairs $\mathop{u_{n}}\limits^{m}(x)$, $\mathop{\lambda_{n}}\limits^{m}$ with the indices $n=1,2,...8$. Vertical axis $\ln \left(\Delta _{n}^{FD} (m)\right)$, horizontal axis $m$.} 
\end{center}
\end{figure}

This test example \eqref{GrindEQ__1_}, \eqref{GrindEQ__2_}, \eqref{GrindEQ__67_} was considered in \cite{58_AttiliLesnic2006,25_SyamSiyyam2009,HPM2010Ex1,HAM2011Ex1,59_Chanane2010} in which the following methods were applied:
Adomian decomposition method (\textbf{ADM}) \cite{58_AttiliLesnic2006}, variational iteration method (\textbf{VIM}) \cite{25_SyamSiyyam2009}, homotopy perturbation method (\textbf{HPM}) \cite{HPM2010Ex1}
homotopy analysis method (\textbf{HAM}) \cite{HAM2011Ex1}, extended sampling method (\textbf{ESM}) \cite{59_Chanane2010}.
Let us compare the behavior of the absolute errors $\Delta _{n}^{FD} (m)$ of FD-method for $m=10,15,20$ with the behavior of the absolute errors $\Delta _{n}^{ADM}$, $\Delta _{n}^{VIM} $, $\Delta _{n}^{HPM} $, $\Delta _{n}^{HAM} $, $\Delta _{n}^{ESM} $  of the corresponding numerical methods from \cite{58_AttiliLesnic2006,25_SyamSiyyam2009,HPM2010Ex1,HAM2011Ex1,59_Chanane2010} as the index number $n$ of eigenvalue increases. 
In Table~2 %\ref{table_MakRom_2} 
 these absolute errors are given. They were calculated using the exact eigenvalues from Table~1 %\ref{table_MakRom_1} 
 and the approximations from \cite{58_AttiliLesnic2006,25_SyamSiyyam2009,HPM2010Ex1,HAM2011Ex1,59_Chanane2010}.

\begin{table}\label{table_MakRom_2}
\begin{center}
{\footnotesize{
\tabcolsep=0.11cm
\begin{tabular}{p{0.1in} p{0.495in} p{0.495in} p{0.495in} p{0.495in} p{0.495in} p{0.495in} p{0.495in} p{0.47in}} \hline 
$n$ & $\Delta _{n}^{FD}(10)$ & $\Delta _{n}^{FD}(15)$ & $\Delta _{n}^{FD}(20)$ & $\Delta _{n}^{ADM}$  & $\Delta _{n}^{VIM}$  & $\Delta _{n}^{HPM}$ & $\Delta _{n}^{HAM}$ & $\Delta _{n}^{ESM}$\\ \hline 
{\footnotesize {(1)}} & {\footnotesize {(2)}} & {\footnotesize {(3)}} & {\footnotesize {(4)}} & {\footnotesize {(5)}} & {\footnotesize {(6)}} & {\footnotesize {(7)}} & {\footnotesize {(8)}} & {\footnotesize {(9)}} \\ \hline 
1 & 4.5e-12 & 9.0e-17 & 1.8e-21 & 4.6e-16 & 5.8e-16 & 1.6e-16  & 4.4e-9  & 2.8e-13\\ 
2 & 4.5e-12 & 9.0e-17 & 1.8e-21 & 4.9e-14 & 1.5e-13 & 1.5e-15 & 1.1e-6 & 3.0e-12 \\ 
3 & 3.4e-16 & 2.3e-22 & 2.4e-28 & 6.5e-13 & 6.3e-13 & 8.5e-13 & 4.4e-5 & 4.2e-11 \\ 
4 & 1.1e-17 & 4.1e-25 & 9.0e-33 & 2.4e-11 & 2.4e-11 & 4.5e-11 & 2.6e-3 & 4.1e-7 \\ 
5 & 3.1e-19 & 7.2e-28 & 1.4e-37 & 7.5e-8 & 7.5e-8  & 2.6e-9 & 0.35 & \\ 
6 & 1.3e-20 & 2.1e-30 & 2.2e-39 & 1.2e-3 & 1.2e-3  & 1.8e-7 & 11.4 & \\ 
7 & 7.8e-22 & 6.4e-33 & 9.2e-42 & 12.7 &  &  &  & \\ 
8 & 6.6e-23 & 8.3e-34 & 5.9e-44 & 215.5 &  &  & &  \\ \hline 
\end{tabular}
\caption{Absolute errors for the first eight eigenvalues $\lambda _{n}$, $n=1,2,...,8$ using the following methods for Example~\ref{example_MakRom_1} :
{columns (2), (3), (4)}  --- FD-method for rank $m=10,15,20$ with error $\Delta _{n}^{FD} (m)$;
{column (5)} --- ADM with error $\Delta _{n}^{ADM}$ \cite{58_AttiliLesnic2006};
{column (6)} --- VIM with error $\Delta _{n}^{VIM} $ \cite{25_SyamSiyyam2009};
{column (7)} --- HPM with error $\Delta _{n}^{HPM} $ \cite{HPM2010Ex1};
{column (8)} --- HAM with error $\Delta _{n}^{HAM} $ \cite{HAM2011Ex1};
{column (9)} --- ESM with error $\Delta _{n}^{ESM} $ \cite{59_Chanane2010}. 
}
}}
\end{center}
\end{table}

One can observe that the convergence rate of each method ADM, VIM, HPM, HAM and ESM rapidly decreases when the eigenvalue index $n$ increases (see columns (5)--(9) in Table~2).
%\ref{table_MakRom_2}). 
The absolute errors $\Delta _{n}^{ADM} $ of the method ADM (variant from \cite{58_AttiliLesnic2006}) for the 7th and 8th  eigenvalues $\lambda _{n}$ with $n=7, 8$ are respectively equal to $12.7$ and $215.5$.
Moreover the convergence of the FD-method increases together with the index $n$ (see columns (2)--(4) in Table~2) and FD-method converges exponentially with respect to rank $m$ (see Figure~1).
The convergence rate is doubled with increase in the rank $m$ (from $10$ to $20$), for example, from $10^{-13}$ (with $m=10$) to $10^{-22}$ (with $m=20$) for $n=1,2$, and from $10^{-24}$ (with $m=10$) to $10^{-45}$ (with $m=20$) for $n=8$.
\end{example}

\newpage
\begin{example}\label{example_MakRom_2}
We consider the eigenvalue problem \eqref{GrindEQ__1_}, \eqref{GrindEQ__2_} with $X=1$, and with potential coefficients in \eqref{GrindEQ__3_} equal to
\begin{equation} \label{GrindEQ__74_} 
q_{0} (x)=x,\;\;\;q_{1} (x)=0,\;\;\;q_{2} (x)=0.
\end{equation}

Developed new symbolic algorithm of the FD-method (see Section~\ref{MakRom_section_6_1}) was applied to compute the approximate solution of the problem \eqref{GrindEQ__1_}, \eqref{GrindEQ__2_}, \eqref{GrindEQ__74_}. The approximate eigenpairs $\mathop{u _{n}}\limits^{m} (x)$, $\mathop{\lambda_{n} }\limits^{m}$ were computed exactly as analytical expressions with respect to the index number $n$ and the input data of the problem under consideration, i.e., we had no rounding errors (see Remark~\ref{remark_MakRom_1}). 
Below we give the eigenvalue corrections $\lambda _{n}^{(j+1)}$ for some first iteration steps $j$ of FD-method:
\begin{equation*}
\begin{split}
&\lambda _{n}^{(0)} =\left(n\pi \right)^{4} ,\;\lambda _{n}^{(1)} =\frac{1}{2} ,\\
&\lambda _{n}^{(2)} =\frac{1}{32\left(n\pi \right)^{4} } -\frac{5}{32\left(n\pi \right)^{6} } +\frac{\cos \left(\pi n\right)-\cosh \left(\pi n\right)}{2\left(n\pi \right)^{7} \sinh \left(\pi n\right)} ,\\
&\lambda _{n}^{(3)} =0,\\
&\lambda _{n}^{(4)} =\frac{11}{20480\left(n\pi \right)^{12} } -\frac{65}{2048\left(n\pi \right)^{14} } -\frac{63e^{2\pi n} -26\cos \left(\pi n\right)e^{\pi n} +63}{128\left(n\pi \right)^{15} \left(e^{2\pi n} -1\right)}\\
 &+\frac{8269e^{2\pi n} +16858\cos \left(\pi n\right)e^{\pi n} +8269}{4096\left(n\pi \right)^{16} \left(e^{2\pi n} +2\cos \left(\pi n\right)e^{\pi n} +1\right)} \\
&-\frac{5\cos \left(\pi n\right)e^{3\pi n} +3e^{2\pi n} -3\cos \left(\pi n\right)e^{\pi n} -5}{16\left(n\pi \right)^{17} \left(\cos \left(\pi n\right)e^{3\pi n} +3e^{2\pi n} +3\cos \left(\pi n\right)e^{\pi n} +1\right)}\\
&-\frac{17\left(e^{2\pi n} -2\cos \left(\pi n\right)e^{\pi n} +1\right)}{32\left(n\pi \right)^{18} \left(e^{2\pi n} +2\cos \left(\pi n\right)e^{\pi n} +1\right)} .
\end{split}
\end{equation*}
Here the command \verb|combine(,trig)| in Maple was used to rewrite the eigenvalue corrections $\lambda _{n}^{(j+1)}$  in a compact form.
The expressions for the eigenfunctions corrections are too cumbersome. Therefore we give the coefficients of \eqref{GrindEQ__15_} only for first two iteration steps of the FD-method with $j=-1,0$: 
\begin{equation*}
\begin{split}
&a_{0}^{(0)}= \sqrt{2} ,\; b_{0}^{(0)} =d_{0}^{(0)} =c_{0}^{(0)} =0,\\
&b_{0}^{(1)} =\frac{\sqrt{2} }{4\left(\pi n\right)^{5} } ,\;b_{1}^{(1)} =\frac{\sqrt{2} }{8\left(\pi n\right)^{3} } ,\;b_{2}^{(1)} =-\frac{\sqrt{2} }{8\left(\pi n\right)^{3} } ,\\
&a_{0}^{(1)} =-\frac{3\sqrt{2} }{16\left(\pi n\right)^{4} } ,\;a_{1}^{(1)} =\frac{3\sqrt{2} }{8\left(\pi n\right)^{4} } ,\;a_{2}^{(1)} =0,\\
&d_{0}^{(1)} =-\frac{\sqrt{2} }{4\left(\pi n\right)^{5} } ,\;c_{0}^{(1)} =-\frac{\sqrt{2} \left(\cos \left(\pi n\right)-\cosh \left(\pi n\right)\right)}{4\left(\pi n\right)^{5} \sinh \left(\pi n\right)} .
\end{split}
\end{equation*}

It should be noted that in this case, the following properties are satisfied for eigenfunction corrections:
\begin{itemize}
\item[--] if the eigenpair index number $n$ is even, then $u_{n}^{(2j+1)} (1-x)=-u_{n}^{(2j+1)} (x)$, $u_{n}^{(2j)} (1-x)=u_{n}^{(2j)} (x)$;
\item[--] if the eigenpair index number $n$ is odd, then $u_{n}^{(2j+1)} (1-x)=u_{n}^{(2j+1)} (x)$, $u_{n}^{(2j)} (1-x)=-u_{n}^{(2j)} (x)$.
\end{itemize}
These properties together with \eqref{GrindEQ__7_}, \eqref{GrindEQ__8_} and Lemma~\ref{MakRom_lemma_1} imply that at the odd-numbered iteration steps $j$ of the FD-method the corrections to the eigenvalues are zero, i. e., $\lambda _{n}^{(2j+1)} =0$, $j=1,2,...$.

In Table~3 the approximations $\mathop{\lambda_{n} }\limits^{10}$ of rank $m=10$ to the exact eigenvalues $\lambda_{n}$ with $n=1,2,3,4,5,10,20,50$ and the norms of the corresponding residuals are given by
\begin{equation} \label{GrindEQ__74_1_} 
\delta _{n} \left(m\right)=\left\| \varphi _{n}^{(m)} (x)\right\| =\left\{\int _{0}^{1}[\varphi 
_{n}^{(m)} (x)]^{2} dx \right\}^{1/2}
\end{equation}
with 
\[\varphi _{n}^{(m)} (x)=\frac{d^{4} \mathop{u_{n}}\limits^{m}  (x)}{dx^{4} } +(x-\mathop{\lambda_n }\limits^{m} )\mathop{u _{n}}\limits^{m} (x).\] 
They are calculated according to the proposed FD-method of rank $m=10$ with the help of the computer algebra system Maple (\verb|Digits=300|).
Figure~2 shows graphs of the approximations $\mathop{u_{n}}\limits^{10}(x)$ to eigenfunctions $u_{n} (x)$ with $n=1,2,3,4,5$.

According to Theorem~\ref{MakRom_theorem_3} the sufficient convergence condition \eqref{GrindEQ__26_} is fulfilled for the eigenpairs with the index $n\ge 2$. For $n = 1$ the FD-method can be divergent. However, as can be seen in Table~4 %\ref{table_MakRom_4} 
and in Figure~3, %\ref{MakRom_fig_3},  
the FD-method converges for $n=1$ too. This means that the conditions of Theorem~\ref{MakRom_theorem_3} can be improved. 
Figure~3 %\ref{MakRom_fig_3} 
shows the 	broken line graphs which were created connecting the data points $\left(m;\ln \left(\delta _{n}(m)\right)\right)$ by lines (see notation \eqref{GrindEQ__74_1_}). 
Figure~3 %\ref{MakRom_fig_3} 
and Table~4 %\ref{table_MakRom_4}
 illustrate the behaviour of the norms of the corresponding residuals $\delta _{n}(m)$ with respect to the rank of FD-method $m$ $(m=1,2,...,10)$ for the indices $n=1,2,3,4,5,10,20,50$, i.e., they illustrate the exponential convergence of the proposed approach for the problem \eqref{GrindEQ__1_}, \eqref{GrindEQ__2_}, \eqref{GrindEQ__74_}.
One can observe that the convergence rate of our method increases together with the index $n$ of the eigenpair $\mathop{u_{n}}\limits^{m}(x)$, $\mathop{\lambda_{n}}\limits^{m}$.

\begin{table}\label{table_MakRom_3}
\begin{center}
{\footnotesize{
\tabcolsep=0.11cm
\begin{tabular}{p{0.1in} p{3.0in} p{0.4in} } \hline 
$n$ & $\mathop{\lambda_{n} }\limits^{10}  $ & $\delta _{n} \left(10\right)$ \\ \hline 
1 &  97.909068819798261176982167541814171360744557739731 & 2.8e-39 \\ 
2 &  1559.0454727668153673091467219850174149875744757492 & 2.7e-39 \\ 
3 &  7890.6363774161879395796364538735759051460151613079 & 5.2e-47 \\ 
4 &  24937.227305908012476430116122759666611086396740215 & 1.2e-51 \\ 
5 &  60881.181896752301770586048651001959246548072513122 & 2.3e-55 \\ 
10 & 974091.41034005627447903500461139135226012366552765 & 1.3e-64 \\ 
20 & 15585455.065440391960236322157494109780226364952145 & 8.4e-74 \\ 
50 & 608806819.46251523277907137706314034909324527027422 & 8.8e-86 \\ \hline 
\end{tabular}
\caption{Approximations $\mathop{\lambda_{n} }\limits^{10} $  to eigenvalues  $\lambda_{n}$ with $n=1,2,3,4,5,10,20,50$ calculated according to the FD-method of ranks $m=10$  and the values of the norms of the corresponding residuals $\delta _{n} \left(m\right)$  for Example~\ref{example_MakRom_2} (see notation \eqref{GrindEQ__74_1_}). 
}
}}
\end{center}
\end{table}

In \cite{RATTANA2013144} matrix methods were developed to approximate the eigenvalues of a fourth order Sturm–-Liouville problem with a kind of fixed boundary conditions. Numerical results for the problem \eqref{GrindEQ__1_}, \eqref{GrindEQ__2_}, \eqref{GrindEQ__74_} were illustrated by using matrix methods such as finite difference method (\textbf{FDM}), modified Numerov's method (\textbf{MNM}), boundary value methods (\textbf{BVM})s of order $p = 6, 8, 10$,  matrix methods FDM*, MNM* and BVMs* of order $p = 6, 8, 10$ with the correction terms (methods denoted with *), as well as ADM and the code SLEUTH (see references in~\cite{RATTANA2013144}). In Table~5 %\ref{table_MakRom_5} 
we illustrate the absolute differences of numerical eigenvalues from \cite[Table~5]{RATTANA2013144} compared with  approximation to the eigenvalues $\mathop{\lambda_{n}}\limits^{10}$, $n=1,2,3,4,5,10,20,50$ which are calculated according to the FD-method of rank $m=10$ and listed in Table~3. %\ref{table_MakRom_3} .
One can observe that the convergence rate of each method FDM, MNM, BVMs of order $p = 6, 8, 10$, ADM and code SLEUTH decreases when the eigenvalue index  $n$ increases (see Table~5), and the convergence rate of each method FDM*, MNM* and BVMs* of order $p = 6, 8, 10$ does not increase, unlike rapid increase of the accuracy of the FD-method with the increasing of the eigenvalue index $n$  (see Table~4 and Figure~3).

\begin{figure}[!t] %float with two figures
\centering
\includegraphics[width=80mm]{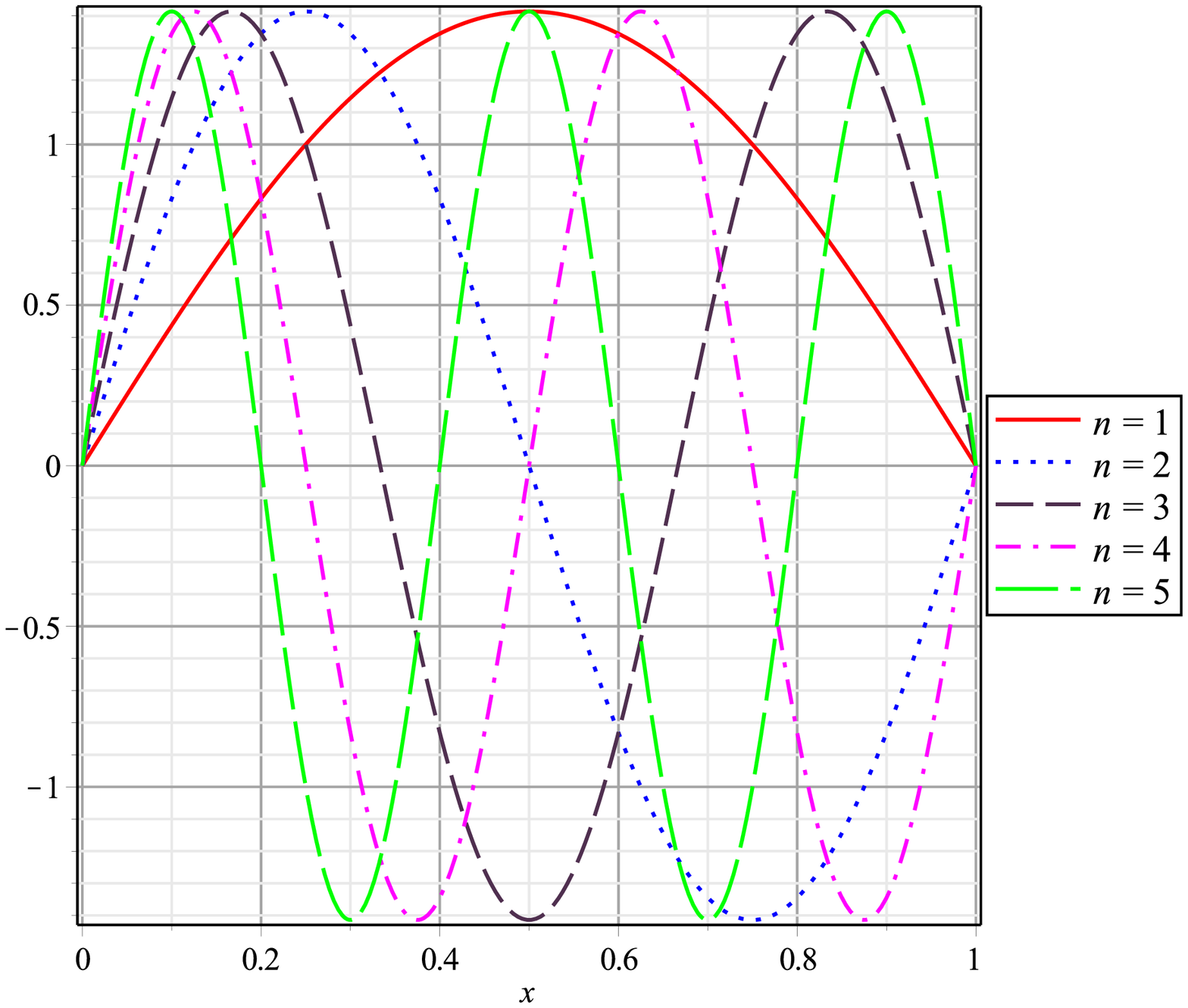}
\caption{Graphs of approximations $\mathop{u_{n}}\limits^{10}(x)$ to first five eigenfunctions $u_{n} (x)$ with $n=1,2,3,4,5$ calculated according to the FD-method of rank $m=10$  for Example~\ref{example_MakRom_2}. Vertical axis $\mathop{u_{n}}\limits^{10}(x)$, horizontal axis $x$.}\bigskip %\label{MakRom_fig_2}
\includegraphics[width=80mm]{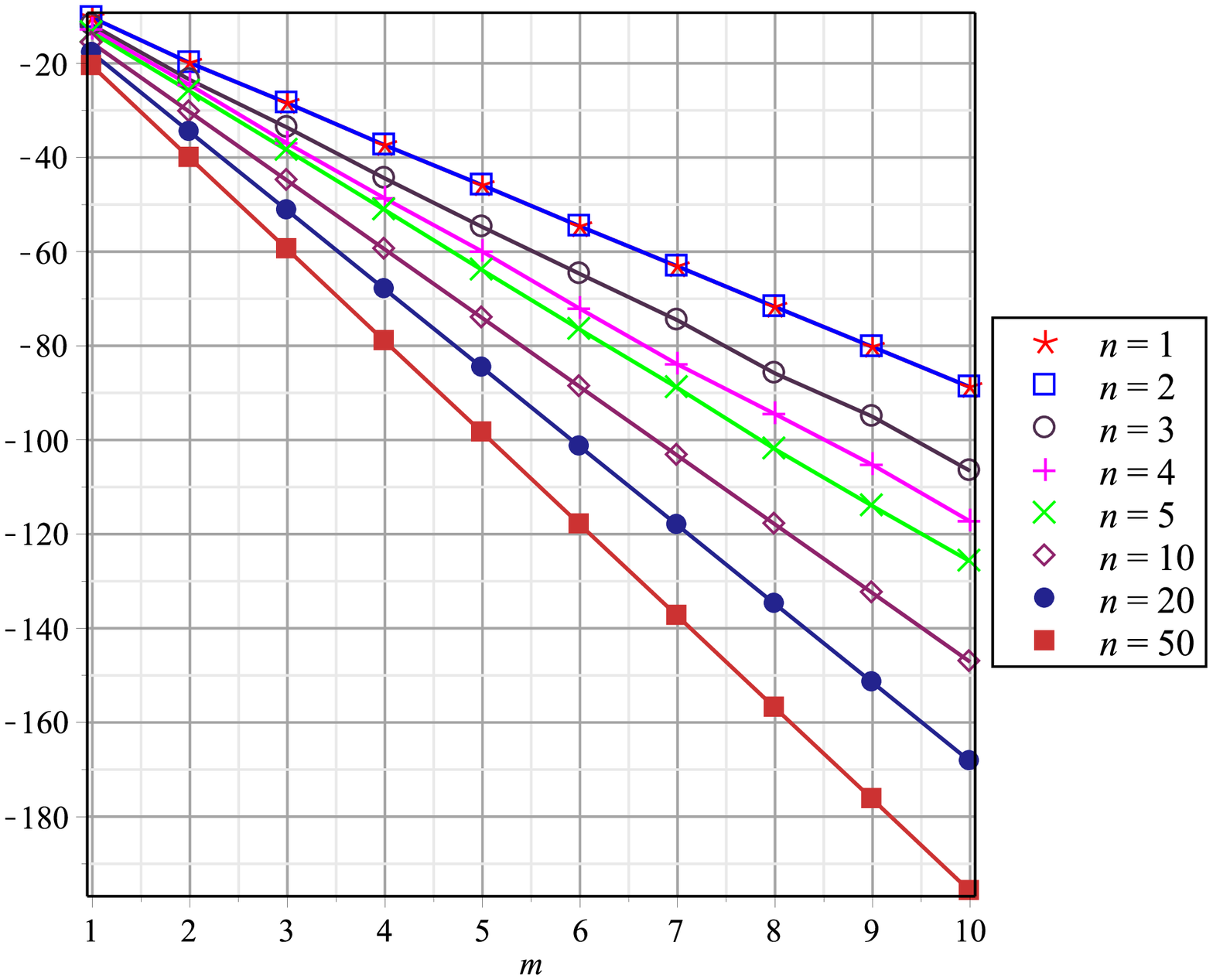}
\caption{Graphs for broken lines joining the data points $\left(m;\ln \left(\delta _{n}(m)\right)\right)$ for Example~\ref{example_MakRom_2}. These graphs illustrates the behaviour of the norms of the corresponding residuals $\delta _{n}(m)$ (see notation \eqref{GrindEQ__74_1_}), i.e., they illustrate the exponential convergence of the FD-method with respect to the rank of FD-method $m$ for the eigenpairs $\mathop{u_{n}}\limits^{m}(x)$, $\mathop{\lambda_{n}}\limits^{m}$ with the indices $n=1,2,3,4,5,10,20,50$.  Vertical axis $\ln \left(\delta _{n}(m)\right)$, horizontal axis $m$.
 }%\label{MakRom_fig_3}
\end{figure}

\newpage

\begin{table}[!t]\label{table_MakRom_4}
\begin{center}
{\footnotesize{
\tabcolsep=0.11cm
\begin{tabular}{p{0.2in} p{0.43in} p{0.43in} p{0.43in} p{0.43in} p{0.43in} p{0.43in} p{0.43in} p{0.43in} } \hline 
$m$  & $\delta _{1} \left(m\right)$ & $\delta _{2} \left(m\right)$ & $\delta _{3} \left(m\right)$ & $\delta _{4} \left(m\right)$ & $\delta _{5} \left(m\right)$ & $\delta _{10} \left(m\right)$ & $\delta _{20} \left(m\right)$ & $\delta _{50} \left(m\right)$ \\ \hline 
1 & 3.3e-5  & 1.7e-5 & 6.7e-6 & 2.7e-6 & 1.3e-6 & 1.5e-7 & 1.8e-8 & 1.1e-9 \\ 
2 & 2.2e-9 & 2.1e-9 & 6.2e-11 & 2.0e-11 & 5.3e-12 & 6.5e-14 & 9.3e-16 & 3.7e-18 \\ 
3 & 4.4e-13  & 2.7e-13 & 2.4e-15 & 8.5e-17 & 1.9e-17 & 3.0e-20 & 5.4e-23 & 1.4e-26 \\  
4 & 5.7e-17  & 5.5e-17 & 4.9e-20 & 7.0e-22 & 5.9e-23 & 1.4e-26 & 2.9e-30 & 4.6e-35 \\  
5 & 1.2e-20 & 6.7e-21 & 1.6e-24 & 8.8e-27 & 1.7e-28 & 6.4e-33 & 1.7e-37 & 1.7e-43 \\  
6 & 1.9e-24  & 1.8e-24 & 7.6e-29 & 4.7e-32 & 5.6e-34 & 2.9e-39 & 8.9e-45 & 5.8e-52 \\  
7 & 3.8e-28 & 2.3e-28 & 4.1e-33 & 3.5e-37 & 2.3e-39 & 1.3e-45 & 5.1e-52 & 2.1e-60 \\  
8 & 7.0e-32 & 6.8e-32 & 5.1e-38 & 9.2e-42 & 5.2e-45 & 6.0e-52 & 2.7e-59 & 7.1e-69 \\  
9 & 1.4e-35 & 8.6e-36 & 5.3e-42 & 1.9e-46 & 3.0e-50 & 2.8e-58 & 1.6e-66 & 2.6e-77 \\  
10 & 2.8e-39 & 2.7e-39 & 5.2e-47 & 1.2e-51 & 2.3e-55 & 1.3e-64 & 8.4e-74 & 8.8e-86 \\ \hline 
\end{tabular}
\caption{Norms of the residuals $\delta _{n}(m)$ (see notation \eqref{GrindEQ__74_1_}) of FD-method of rank $m=1,2,...,10$ for the approximations to eigenpairs $\mathop{u_{n}}\limits^{m}(x)$, $\mathop{\lambda_{n}}\limits^{m}$ with the indices $n=1,2,3,4,5,10,20,50$ for Example~\ref{example_MakRom_2}.
}
}}
\end{center}
\end{table}

$$\;$$

\begin{table}[!t]\label{table_MakRom_5}
\begin{center}
{\footnotesize{
\tabcolsep=0.11cm
\begin{tabular}{p{0.5in} p{0.43in} p{0.43in} p{0.43in} p{0.43in} p{0.43in} p{0.43in} p{0.43in} p{0.43in} } \hline 
Method  & $n=1$ & $n=2$ & $n=3$ & $n=4$ & $n=5$ & $n=10$ & $n=20$ & $n=50$ \\ \hline 
FDM & 0.4e-2 & 2.5e-1 & 2.9 & 1.6+1 & 6.2e+1 & 4.0e+3 & 2.5e+5 & 5.9e+7 \\ 
MNM & 9.9e-7 & 2.9e-6 & 5.3e-5 & 5.3e-4 & 3.2e-3 & 8.1e-1 & 2.1e+2 & 3.6e+5 \\ 
order6 & 5.2e-7 & 5.8e-7 & 1.1e-6 & 4.1e-7 & 4.7e-6 & 4.7e-3 & 4.9 & 4.8e+3 \\ 
order8 & 8.9e-5 & 1.5e-3 & 7.4e-3 & 2.3e-2 & 5.7e-2 & 9.1e-1 & 1.4e+1 & 4.1e+3 \\ 
order10 & 1.9e-5 & 3.0e-4 & 1.5e-3 & 4.7e-3 & 1.2e-2 & 1.8e-1 & 2.9 & 5.4e+2 \\ 
FDM* & 1.8e-6 & 1.2e-6 & 3.6e-6 & 2.0e-6 & 6.8e-7 & 1.7e-6 & 4.2e-6 & 1.2e-6 \\ 
MNM* & 9.1e-7 & 1.1e-5 & 4.4e-7 & 5.5e-6 & 4.5e-7 & 9.1e-7 & 6.7e-6 & 4.0e-6 \\ 
order6* & 3.0e-6 & 8.3e-6 & 3.9e-6 & 3.7e-6 & 2.1e-6 & 5.7e-6 & 7.6e-6 & 4.8e-6 \\ 
order8* & 7.9e-8 & 6.6e-7 & 1.5e-6 & 1.6e-7 & 2.0e-6 & 5.5e-7 & 3.7e-6 & 4.3e-7 \\ 
order10* & 3.1e-6 & 2.1e-6 & 2.5e-6 & 1.1e-6 & 2.4e-6 & 1.9e-6 & 2.9e-6 & 2.7e-7 \\ 
ADM & 1.2e-15 & 3.7e-13 & 2.0e-8 & 1.0e-2 & 2.8+2 &  &  &   \\ 
SLEUTH & 2.0e-8 & 2.8e-6 & 2.6e-6 & 5.9e-6 & 3.2e-6 & 3.4e-4 & 3.5e-2 & 5.4e-1 \\ \hline 
\end{tabular}
\caption{Absolute differences in the approximations to eigenvalues from \cite[Table~5]{RATTANA2013144} with respect to the numerical eigenvalues $\mathop{\lambda_{n}}\limits^{10}$, $n=1,2,3,4,5,10,20,50$  calculated according to the FD-method of rank $m=10$ for Example~\ref{example_MakRom_2}.
}
}}
\end{center}
\end{table}

\end{example}

\newpage
\section{Conclusions}\label{MakRom_section_8}%8
\textbf{Results}. 
In this article a new symbolic algorithmic implementation of the functional-discrete (FD-) method is developed and justified for the fourth order Sturm--Liouville problem (see numerical algorithm from Section~\ref{MakRom_section_6_1}). We consider the eigenvalue problem on a finite interval $[0,X]$ in the Hilbert space $L_{2} \left(0,X\right)$ for the fourth order ordinary differential equation \eqref{GrindEQ__1_} with polynomial coefficients \eqref{GrindEQ__3_} and boundary conditions  \eqref{GrindEQ__2_}.
The sufficient conditions of an exponential convergence rate of FD-method are received (see Theorem~\ref{MakRom_theorem_3}). The obtained estimates of the absolute errors of FD-method \eqref{GrindEQ__27_}, \eqref{GrindEQ__28_}  significantly improve the accuracy of the estimates obtained earlier in \cite{GMR-GavrMakPop2010}.
The theoretical results are illustrated by numerical examples \ref{example_MakRom_1} and \ref{example_MakRom_2} in which the numerical results obtained with the FD-method are compared with the numerical test results obtained with other existing numerical techniques \cite{58_AttiliLesnic2006,25_SyamSiyyam2009,HPM2010Ex1,HAM2011Ex1,59_Chanane2010,RATTANA2013144}. 

\textbf{Features of implementation}. 
The obtained algorithm is symbolic and operates with the decomposition coefficients \eqref{GrindEQ__66_} 
of the eigenfunction corrections $u_{n}^{(j+1)} (x)$ in some basis on interval $\left[0,X\right]$ (see Lemma~\ref{MakRom_lemma_2}). Unlike the symbolic algorithm from \cite{Conf_MakRom2015} and traditional algorithm from \cite{GMR-GavrMakPop2010,GavrMakRom20152017}, presented approach produces explicit recursive formulas for the coefficients in (\ref{GrindEQ__15_}) which are corresponding elements of the column vectors \eqref{GrindEQ__65_1_} and \eqref{GrindEQ__65_2_}. These coefficients are represented recursively through the coefficients and quantities computed at previous steps of FD-method.

\textbf{Unique advantages}.
Proposed symbolic algorithm of the simplest variant of the FD-method for problem \eqref{GrindEQ__1_}--\eqref{GrindEQ__3_} will  always be convergent beginning with some eigenvalue index number $n_{0} $ (it can be large enough) with estimates of the absolute errors \eqref{GrindEQ__27_} and \eqref{GrindEQ__28_} (see Theorem~\ref{MakRom_theorem_3}). %\ref{MakRom_fig_3}).
 This means that by using the simplest variant of the FD-method one can obtain the asymptotic formulas for eigenvalues and eigenfunctions. 

The approximate eigenpairs $\mathop{u _{n}}\limits^{m} (x)$, $\mathop{\lambda_{n} }\limits^{m}$ are computed exactly as analytical expressions and there are no rounding errors (see Remark~\ref{remark_MakRom_1}). Proposed symbolic algorithm uses only the algebraic operations and basic operations on $(2\times 1)$ column vectors and $(2\times 2)$ matrices. Presented method does not require solving any boundary value problems and computations of any integrals, unlike the previous variants of FD-method from \cite{GMR-GavrMakPop2010,GavrMakRom20152017}. Substituting into the obtained analytical expressions the given value $n$ and the numerical values of input data, we find numerical values of the corresponding approximations $\mathop{\lambda_{n} }\limits^{m}$, $\mathop{u_{n}}\limits^{m}  (x)$. 

 FD-method converges exponentially with respect to rank $m$. Moreover the convergence of the FD-method increases together with the index $n$, unlike the accuracy degradation of other existing numerical techniques with the increasing of the eigenvalue index $n$.
 
%\section*{Acknowledgement}%9

%\newpage
\section*{Appendix A}\label{MakRom_append_A}
Analytical formulas which are used for computation in the proposed numerical algorithm (see the step 4 and 11 of the numerical algorithm from Section~\ref{MakRom_section_6_1}): 
\begin{equation*}
\begin{split}
&F_{n,\cos }^{(j+1)} (x)=\sum _{t=0}^{M(j)}x^{t}  \sum _{s=\left. \left. \right]\kern-0.15em\right]\frac{t}{r+1} \left[\kern-0.15em\left[\right. \right. }^{j}\lambda_n ^{(j+1-s)}  b_{n,t}^{(s)} \\
&+\sum _{t=0}^{M(j+1)-1}x^{t}  \sum _{l=\max (0,t-M(j))}^{\min (r,t)}\left(b_{n,t-l}^{(j)} \left(-A_{l} +C_{l} \left(\frac{\pi n}{X} \right)^{2} \right)-a_{n,t-l}^{(j)} B_{l} \frac{\pi n}{X} \right) \\
&-\sum _{t=0}^{M(j+1)-2}x^{t}  \sum _{l=\max (0,t-M(j)+1)}^{\min (r,t)}\left(b_{n,t-l+1}^{(j)} B_{l} +a_{n,t-l+1}^{(j)} C_{l} 2\frac{\pi n}{X} \right)\left(t-l+1\right) \\
&-\sum _{t=0}^{M(j+1)-3}x^{t}  \sum _{l=\max (0,t-M(j)+2)}^{\min (r,t)}b_{n,t-l+2}^{(j)} C_{l} \left(t-l+2\right)\left(t-l+1\right) ,
\end{split}
\end{equation*}
\begin{equation*}
\begin{split}
&F_{n,\sin }^{(j+1)} (x)=\sum _{t=0}^{M(j)}x^{t}  \sum _{s=\left. \left. \right]\kern-0.15em\right]\frac{t}{r+1} \left[\kern-0.15em\left[\right. \right. }^{j}\lambda_n ^{(j+1-s)}  a_{n,t}^{(s)}\\
&+\sum _{t=0}^{M(j+1)-1}x^{t}  \sum _{l=\max (0,t-M(j))}^{\min (r,t)}\left(n,a_{n,t-l}^{(j)} \left(-A_{l} +C_{l} \left(\frac{\pi n}{X} \right)^{2} \right)+b_{n,t-l}^{(j)} B_{l} \frac{\pi n}{X} \right)\\
&-\sum _{t=0}^{M(j+1)-2}x^{t}  \sum _{l=\max (0,t-M(j)+1)}^{\min (r,t)}\left(a_{n,t-l+1}^{(j)} B_{l} -b_{n,t-l+1}^{(j)} C_{l} 2\frac{\pi n}{X} \right)\left(t-l+1\right)\\
&-\sum _{t=0}^{M(j+1)-3}x^{t}  \sum _{l=\max (0,t-M(j)+2)}^{\min (r,t)}a_{n,t-l+2}^{(j)} C_{l} \left(t-l+2\right)\left(t-l+1\right) ,
\end{split}
\end{equation*}
\begin{equation*}
\begin{split}
&F_{n,\cosh }^{(j+1)} (x)=\sum _{t=0}^{M(j-1)}x^{t}  \sum _{s\left. \left. =\right]\kern-0.15em\right]\frac{t}{r+1} \left[\kern-0.15em\left[+\right. \right. 1}^{j}\lambda_n ^{(j+1-s)}  d_{n,t}^{(s)} \\
&-\sum _{t=0}^{M(j)-1}x^{t}  \sum _{l=\max (0,t-M(j-1))}^{\min (r,t)}\left(d_{n,t-l}^{(j)} \left(A_{l} +C_{l} \left(\frac{\pi n}{X} \right)^{2} \right)+c_{n,t-l}^{(j)} B_{l} \frac{\pi n}{X} \right)\\
&-\sum _{t=0}^{M(j)-2}x^{t}  \sum _{l=\max (0,t-M(j-1)+1)}^{\min (r,t)}\left(d_{n,t-l+1}^{(j)} B_{l} +c_{n,t-l+1}^{(j)} C_{l} 2\frac{\pi n}{X} \right)\left(t-l+1\right) \\
&-\sum _{t=0}^{M(j)-3}x^{t}  \sum _{l=\max (0,t-M(j-1)+2)}^{\min (r,t)}d_{n,t-l+2}^{(j)} C_{l} \left(t-l+2\right)\left(t-l+1\right) ,
\end{split}
\end{equation*}
\begin{equation*}
\begin{split}
&F_{n,\sinh }^{(j+1)} (x)=\sum _{t=0}^{M(j-1)}x^{t}  \sum _{s\left. \left. =\right]\kern-0.15em\right]\frac{t}{r+1} \left[\kern-0.15em\left[+\right. \right. 1}^{j}\lambda_n ^{(j+1-s)}  c_{n,t}^{(s)}\\
&-\sum _{t=0}^{M(j)-1}x^{t}  \sum _{l=\max (0,t-M(j-1))}^{\min (r,t)}\left(c_{n,t-l}^{(j)} \left(A_{l} +C_{l} \left(\frac{\pi n}{X} \right)^{2} \right)+d_{n,t-l}^{(j)} B_{l} \frac{\pi n}{X} \right)\\
&-\sum _{t=0}^{M(j)-2}x^{t}  \sum _{l=\max (0,t-M(j-1)+1)}^{\min (r,t)}\left(c_{n,t-l+1}^{(j)} B_{l} +d_{n,t-l+1}^{(j)} C_{l} 2\frac{\pi n}{X} \right)\left(t-l+1\right)\\
&-\sum _{t=0}^{M(j)-3}x^{t}  \sum _{l=\max (0,t-M(j-1)+2)}^{\min (r,t)}c_{n,t-l+2}^{(j)} C_{l} \left(t-l+2\right)\left(t-l+1\right) .
\end{split}
\end{equation*}
These formulas enter into the expression \eqref{GrindEQ__29_1}. Here $\left. \left. \right]\kern-0.15em\right]y\left[\kern-0.15em\left[\right. \right. $ is the smallest integer greater than or equal to a real number $y$ (this is the function \verb|ceil(y)| in Maple).

\section*{Appendix B}\label{MakRom_append_B}
Analytical expressions for the integrals $\alpha _{n,t}$, $\beta _{n,t}$, $\eta _{n,t}$, $\mu _{n,t}$ which are used in \eqref{GrindEQ__66_2_} and exactly calculated in the case $t=0,1,...$  \cite{GradshteynRyzhik2014} (see step 2 in Section~\ref{MakRom_section_6_1}):
\begin{equation*}
\begin{split}
\alpha _{n,t} &=\int _{0}^{X}\xi ^{t} \sin ^{2} \left(\frac{\pi n}{X} \xi \right)d\xi  =\frac{1}{2} \frac{X^{t+1} }{t+1}
-\frac{1}{2} t!X^{t+1} \sum _{k=0}^{t-1} \frac{1}{(t-k)!(2\pi n)^{k+1} } \sin \left(\frac{\pi k}{2} \right),
\end{split}
\end{equation*}
\begin{equation*}
\begin{split}
\beta _{n,t} =\frac{1}{2} \int _{0}^{X}\xi ^{t} \sin \left(\frac{2\pi n}{X} \xi \right)d\xi  =-t!X^{t+1} \frac{1}{2} \sum _{k=0}^{t-1} \frac{1}{(t-k)!(2\pi n)^{k+1} } \cos \left(\frac{\pi k}{2} \right),
\end{split}
\end{equation*}
\begin{equation*}
\begin{split}
\eta _{n,t} &=\int _{0}^{X}\xi ^{t} \sin \left(\frac{\pi n}{X} \xi \right)\cosh \left(\frac{\pi n}{X} \xi \right) d\xi 
=\frac{X^{t+1} t!}{\left(\sqrt{2} \pi n\right)^{t+1} } \cos \left(\frac{\pi t}{2} \right)\cos \left(\frac{\pi \left(t+1\right)}{4} \right)\\
\end{split}
\end{equation*}
\begin{equation*}
\begin{split}
&-\sum _{k=0}^{t}\frac{t!X^{t+1} \cos \left(\pi n\right)}{\left(t-k\right)!\left(\sqrt{2} \pi n\right)^{k+1} }  \left[\cos \left(\frac{\pi \left(k+1\right)}{4} \right)\right. \cos \left(\frac{\pi k}{2} \right)\cosh \left(\pi n\right)\\
\end{split}
\end{equation*}
\begin{equation*}
\begin{split}
&-\left. \sin \left(\frac{\pi \left(k+1\right)}{4} \right)\sin \left(\frac{\pi k}{2} \right)\sinh \left(\pi n\right)\right],
\end{split}
\end{equation*}
\begin{equation*}
\begin{split}
\mu _{n,t} &=\int _{0}^{X}\xi ^{t} \sin \left(\frac{\pi n}{X} \xi \right)\sinh \left(\frac{\pi n}{X} \xi \right) d\xi 
=-\frac{X^{t+1} t!}{\left(\sqrt{2} \pi n\right)^{t+1} } \sin \left(\frac{\pi t}{2} \right)\sin \left(\frac{\pi \left(t+1\right)}{4} \right)\\
&+\sum _{k=0}^{t}\frac{t!X^{t+1} \cos \left(\pi n\right)}{\left(t-k\right)!\left(\sqrt{2} \pi n\right)^{k+1} }  \left[\sin \left(\frac{\pi \left(k+1\right)}{4} \right)\sin \left(\frac{\pi k}{2} \right)\cosh \left(\pi n\right)\right.\\
&\left.-\cos \left(\frac{\pi \left(k+1\right)}{4} \right) \cos \left(\frac{\pi k}{2} \right)\sinh \left(\pi n\right)\right].
\end{split}
\end{equation*}

\section*{Appendix C}\label{MakRom_append_C} 

The formula for the corrections of eigenvalues (see steps 3,10 in Section~\ref{MakRom_section_6_1}):
\begin{equation*}
\begin{split}
\lambda _{n}^{(j+1)} =&a_{n,0}^{(0)} \cdot \left(\sum _{t=0}^{M(j+1)-1} \right.\sum _{l=\max (0,t-M(j))}^{\min (r,t)} \left[\frac{\pi n}{X} B_{l} \left(\beta _{n,t} a_{n,t-l}^{(j)} -\alpha _{n,t} b_{n,t-l}^{(j)} \right)\right.\\
&-\left.\left(\beta _{n,t} b_{n,t-l}^{(j)} +\alpha _{n,t} a_{n,t-l}^{(j)} \right)\left(-A_{l} +\left(\frac{\pi n}{X} \right)^{2} C_{l} \right)\right]+\\
\sum _{t=0}^{M(j+1)-2} &\sum _{l=\max (0,t-M(j)+1)}^{\min (r,t)} \left(t-l+1\right)\left[B_{l} \left(\beta _{n,t} b_{n,t-l+1}^{(j)} +\alpha _{n,t} a_{n,t-l+1}^{(j)} \right)\right.\\
&\left.+2\frac{\pi n}{X} C_{l} \left(\beta _{n,t} a_{n,t-l+1}^{(j)} -\alpha _{n,t} b_{n,t-l+1}^{(j)} \right)\right]\\
+\sum _{t=0}^{M(j+1)-3} &\sum _{l=\max (0,t-M(j)+2)}^{\min (r,t)} \left(t-l+1\right)\left(t-l+2\right)C_{l} \left(\beta _{n,t} b_{n,t-l+2}^{(j)} +\alpha _{n,t} a_{n,t-l+2}^{(j)} \right)\\
+\sum _{t=0}^{M(j)-1} &\sum _{l=\max (0,t-M(j-1))}^{\min (r,t)} \left[\frac{\pi n}{X} B_{l} \left(\eta _{n,t} c_{n,t-l}^{(j)} +\mu _{n,t} d_{n,t-l}^{(j)} \right)\right.\\
&\left.+\left(\eta _{n,t} d_{n,t-l}^{(j)} +\mu _{n,t} c_{n,t-l}^{(j)} \right)\left(A_{l} +\left(\frac{\pi n}{X} \right)^{2} C_{l} \right)\right]\\
+\sum _{t=0}^{M(j)-2}  &\sum _{l=\max (0,t-M(j-1)+1)}^{\min (r,t)} \left(t-l+1\right)\left[B_{l} \left(\eta _{n,t} d_{n,t-l+1}^{(j)} +\mu _{n,t} c_{n,t-l+1}^{(j)} \right)\right.\\
&\left.+2\frac{\pi n}{X} C_{l} \left(\eta _{n,t} c_{n,t-l+1}^{(j)} +\mu _{n,t} d_{n,t-l+1}^{(j)} \right)\right]\\
\end{split}
\end{equation*}
\begin{equation*}
\begin{split}
+\sum _{t=0}^{M(j)-3} &\left.\sum _{l=\max (0,t-M(j-1)+2)}^{\min (r,t)} \left(t-l+1\right)\left(t-l+2\right)C_{l} \left(\eta _{n,t} d_{n,t-l+2}^{(j)} +\mu _{n,t} c_{n,t-l+2}^{(j)} \right)\right).
\end{split}
\end{equation*}

\newpage
\section*{References}

{\footnotesize{
\bibliography{mybibfile}
}}

\end{document}